\documentclass[reqno]{amsart}
\usepackage{amssymb, latexsym, euscript}

\newtheorem{theorem}{Theorem}[section]
\newtheorem{proposition}[theorem]{Proposition}
\newtheorem{cor}[theorem]{Corollary}
\newtheorem{lemma}[theorem]{Lemma}
\newtheorem{definition}[theorem]{Definition}
\newtheorem{corollary}[theorem]{Corollary}

\theoremstyle{definition}
\newtheorem{example}[theorem]{Example}
\newtheorem{remark}[theorem]{Remark}

\numberwithin{equation}{section}

   \def\sH{{\mathfrak H}}

   \def\sT{{\mathfrak T}}   \def\sU{{\mathfrak U}}

\def\dA{{\mathbb A}}      \def\dC{{\mathbb C}}

      \def\dR{{\mathbb R}}
      \def\dU{{\mathbb U}}

\def\cD{{\mathcal D}}      
   \def\cH{{\mathcal H}}

      \def\cR{{\mathcal R}}

\def\cY{{\mathcal Y}}

\def\wt#1{{{\widetilde #1} }}

\def\wh#1{{{\widehat #1} }}

\def\bm\chi{\mbox{\boldmath$\chi$}}

\def\ker{{\rm ker\,}}
\def\ran{{\rm ran\,}}

\def\dim{{\rm dim\,}}

\let\xker=\ker \def\ker{{\xker\,}}

\def\cmr{{\dC \setminus \dR}}

\begin{document}

\title{On Symmetries in the Theory of Finite Rank
Singular Perturbations}
\author[S.~Hassi]{Seppo~Hassi}
\author[S.~Kuzhel]{Sergii~Kuzhel}

\address{Department of Mathematics and Statistics \\
University of Vaasa \\
P.O. Box 700, 65101 Vaasa \\
Finland} \email{sha@uwasa.fi}

\address{Institute of Mathematics of the National
Academy of Sciences of Ukraine \\
3 Tere\-shchenkivska Street, 01601, Kiev-4 \\
Ukraine} \email{kuzhel@imath.kiev.ua}



\begin{abstract}
For a nonnegative self-adjoint operator $A_0$ acting on a Hilbert
space $\mathfrak{H}$ singular perturbations of the form $A_0+V, \
V=\sum_{1}^{n}{b}_{ij}<\psi_j,\cdot>\psi_i$ are studied under some
additional requirements of symmetry imposed on the initial operator
$A_0$ and the singular elements $\psi_j$. A concept of symmetry is
defined by means of a one-parameter family of unitary operators
$\sU$ that is motivated by results due to R. S. Phillips. The
abstract framework to study singular perturbations with symmetries
developed in the paper allows one to incorporate physically
meaningful connections between singular potentials $V$ and the
corresponding self-adjoint realizations of $A_0+V$. The results are
applied for the investigation of singular perturbations of the
Schr\"{o}dinger operator in $L_2(\dR^3)$ and for the study of a
(fractional) \textsf{p}-adic Schr\"{o}dinger type operator with
point interactions.
\end{abstract}

\maketitle
\section{Introduction}
Let $A_0$ be an unbounded nonnegative self-adjoint operator acting
on a Hilbert space $\mathfrak{H}$ and let $
 \mathfrak{H}_2(A_0)\subset\mathfrak{H}_1(A_0)\subset\mathfrak{H}
 \subset\mathfrak{H}_{-1}(A_0)\subset\mathfrak{H}_{-2}(A_0)$
 be the standard scale of Hilbert spaces associated with $A_0$.
More precisely,
\begin{equation}\label{ee1}
 {\mathfrak H}_k(A_0)=\mathcal{D}(A_0^{k/2}), \quad k=1,2,
\end{equation}
equipped with the norm $\|u\|_k=\|(A_0+I)^{k/2}u\|$. The dual spaces
${\mathfrak H}_{-k}(A_0)$ can be defined as the completions of
${\mathfrak H}$ with respect to the norms
$\|u\|_{-k}=\|(A_0+I)^{-k/2}u\|$ \ $(u \in\mathfrak{H})$. The
resolvent operator $(A_0+I)^{-1}$ can be continuously extended to an
isometric mapping $(\mathbb{A}_0+I)^{-1}$ from ${\mathfrak
H}_{-2}(A_0)$ onto ${\mathfrak H}$ and the relation
\begin{equation}\label{ada5}
 <\psi,u>=((A_0+I)u,(\mathbb{A}_0+I)^{-1}\psi),
 \quad u\in{\mathfrak H}_{2}(A_0)
\end{equation}
enables one to identify the elements $\psi\in{\mathfrak
H}_{-2}(A_0)$ as linear functionals on ${\mathfrak H}_{2}(A_0)$.

Consider the heuristic expression
\begin{equation}\label{ne3}
 A_0 +\sum_{i,j=1}^{n}{b}_{ij}<\psi_j,\cdot>\psi_i, \quad
 b_{ij}\in\mathbb{C}, \quad n\in\mathbb{N},
\end{equation}
where elements $\psi_j$ $(1\leq{j}\leq{n})$ form a linearly
independent system in ${\mathfrak  H}_{-2}(A_0)$. In what follows it
is supposed that the linear span $\mathcal{X}$ of
$\{\psi_j\}_{j=1}^n$  satisfies the condition
$\mathcal{X}\cap{\mathfrak H}=\{0\}$, i.e., elements $\psi_j$ are
${\mathfrak H}$-independent. In this case, the perturbation
 $V=\sum_{i,j=1}^{n}{b}_{ij}<\psi_j,\cdot>\psi_i$ is said to be singular
and the formula
\begin{equation}\label{e7}
 A_{\mathrm{sym}}={A}_0\upharpoonright{{\mathcal{D}}(A_{\mathrm{sym}})},
  \quad {\mathcal{D}}(A_{\mathrm{sym}})
    =\{\,u\in{{\mathcal{D}}(A_0)}  :\,  <\psi_j, u>=0, \ 1\leq{j}\leq{n} \,\}
\end{equation}
determines a closed densely defined symmetric operator in
${\mathfrak H}$.

In the theory of singular perturbations, cf. e.g. \cite{AL, AL1,
HS}, each intermediate extension $A$ of $A_{\mathrm{sym}}$, i.e.,
$A_{\mathrm{sym}}\subset{A}\subset{A}_{\mathrm{sym}}^*$, can be
viewed to be singularly perturbed with respect to $A_0$ and, in
general, such an extension can be regarded as an
operator-realization of (\ref{ne3}) in $\mathfrak{H}$. In this
context, the natural question arises whether and how one could
establish a physically meaningful correspondence between the
parameters $b_{ij}$ of the singular potential $V$ and the
intermediate extensions of $A_{\mathrm{sym}}$. The investigation of
this problem is one of goals of the present paper. In the approach
developed in \cite{AL1, A4} one considers an operator realization
$A$ of (\ref{ne3}) by setting
\begin{equation}\label{lesia40}
 A=\mathbb{A}_\mathbf{R}\upharpoonright{\mathcal{D}(A)}, \quad
 \mathcal{D}(A)=\{\,f\in\mathcal{D}(A_{\mathrm{sym}}^*) :\,
 \mathbb{A}_\mathbf{R}f\in\mathfrak{H}\,\},
\end{equation}
where
\begin{equation}\label{tat44}
 \mathbb{A}_\mathbf{R}=
 \mathbb{A}_0+\sum_{i,j=1}^n{b}_{ij}<\psi_j^{\mathrm{ex}},\cdot>\psi_i
\end{equation}
is seen as a regularization of (\ref{ne3}).

Formula (\ref{tat44}) involves a construction of the extended
functionals $<\psi_j^{\mathrm{ex}}, \cdot>$ defined on
$\mathcal{D}(A_{\mathrm{sym}}^*)$. These functionals are uniquely
determined by the choice of a Hermitian matrix
$\mathbf{R}=(r_{jp})_{j,p=1}^n$. Since for elements
$\psi\in\mathcal{X}\cap\mathfrak{H}_{-1}(A_0)$ the functionals
$<\psi, \cdot>$ admit extensions by continuity onto
$\mathfrak{H}_{1}(A_0)\cap\mathcal{D}(A_{\mathrm{sym}}^*)$, a lot of
natural restrictions appears in the choice of $\mathbf{R}$. For
their preservation the concept of admissible matrices $\mathbf{R}$
for the regularization of (\ref{ne3}) has been introduced in
\cite[Definition 3.1.2]{A4}. However, this definition involves
certain spectral measures and, in what follows, their calculation
will be avoided. In fact, an equivalent concept of
\textit{admissible operators} is introduced in the form convenient
for the further studies in the present paper.

If the singular potential $V$ in (\ref{ne3}) is not form-bounded
(i.e., $\mathcal{X}\not\subset{{\mathfrak H}_{-1}(A_0)}$), then an
admissible operator cannot be determined uniquely and one needs to
impose some extra assumptions to achieve the uniqueness. For
instance, in many applications, the condition of extremality
\cite{Ar, AHSS} allows one to select a unique admissible operator
(see Theorem \ref{did3}). It should be noted that the concept of
extremality is physically reasonable. For example, extremal
operators determine free evolutions in the Lax--Phillips scattering
theory \cite{Ku5}.

Another approach inspired by \cite{AL1, A4, KP1} deals with the
preservation of initially existing symmetries of singular elements
$\psi_j$ in the definition of the extended functionals
$\psi_j^{\mathrm{ex}}$. To study this problem in an abstract
framework, one needs to define the notion of symmetry for the
unperturbed operator $A_0$ and for the singular elements $\psi_j$ in
(\ref{ne3}). Generalizing the ideas suggested in \cite{AL1, KO, PH},
the required definitions will be formulated here as follows:

Let ${\mathfrak{T}}$ be a subset of the real line $\mathbb{R}$ and
let ${\mathfrak U}=\{U_t\}_{t\in{\mathfrak{T}}}$ be a one-parameter
family of unitary operators acting on ${\mathfrak H}$ with the
following property:
\begin{equation}\label{e1}
 U_t\in{\mathfrak U} \iff U^*_t\in{\mathfrak U}
\end{equation}

\begin{definition}{\cite{MS12}}\label{dad1}
A linear operator $A(\not=0)$ acting in ${\mathfrak H}$ is said to
be $p(t)$-homogeneous with respect to ${\mathfrak U}$ if there
exists a real function $p(t)$ defined on ${\mathfrak{T}}$ such that
 \begin{equation}\label{ee6}
 U_tA=p(t)AU_t, \quad \forall{t}\in\mathfrak{T}.
 \end{equation}
\end{definition}

In other words, the set ${\mathfrak U}$ determines the structure of a symmetry
and the property of $A$ to be $p(t)$-homogeneous with respect to ${\mathfrak U}$
means that $A$ possesses a certain symmetry with respect to ${\mathfrak U}$.

\begin{definition}{\cite{MS12}}\label{dad2}
A singular element $\psi\in{\mathfrak
H}_{-2}(A_0)\setminus{\mathfrak  H}$ is said to be
$\xi(t)$-invariant with respect to ${\mathfrak U}$ if there exists a
real function $\xi(t)$ defined on ${\mathfrak{T}}$  such that
\begin{equation}\label{tat10}
 {\mathbb{U}}_t\psi=\xi(t)\psi, \quad \forall{t}\in\mathfrak{T},
\end{equation}
where ${\mathbb{U}}_t$ is the continuation of $U_t$ onto ${\mathfrak
H}_{-2}(A_0)$  (see Section~\ref{sec4} for details).
\end{definition}

The main aim of the paper is to study (\ref{ne3}) assuming that the
initial operator $A_0$ is $p(t)$-homogeneous and the singular
elements $\psi_j$ are $\xi_j(t)$-invariant with respect to
${\mathfrak U}$. It appears that the preservation of
$\xi_j(t)$-invariance for the extended functionals
$<\psi_j^{\mathrm{ex}}, \cdot>$  is equivalent to the
$p(t)$-homogeneity of the operator $\widetilde{A}$ which is used for
the regularization of (\ref{ne3}) (Theorem~\ref{t5}). Combining this
result with the complete description of admissible operators
(Theorem \ref{did121}) allows one to select a unique admissible
operator by imposing the condition of $p(t)$-homogeneity (Theorems
\ref{new2007}, \ref{did11}). One of interesting properties
discovered here is the possibility to get the Friedrichs and the
Krein-von Neumann extension (and more generally, all
$p(t)$-homogeneous self-adjoint extensions transversal to $A_0$) as
solutions of a system of equations involving the functions $p(t)$
and $\xi_j(t)$ (Corollary \ref{new28}, Proposition \ref{pip1}).

The choice of a $p(t)$-homogeneous admissible operator for the
regularization of (\ref{ne3}) immediately gives a new specific
relation for the corresponding Weyl function $\mathbf{M}(z)$
(Theorem \ref{pop1}) and enables one to establish simple relations
involving the functions $p(t)$ and $\xi_j(t)$, and the properties of
operator realizations of (\ref{ne3}) (Theorem \ref{p5}, Proposition
\ref{new2008}).

It is well known, see e.g. \cite{ADK, CF, KP, KP1} that the
Schr\"{o}dinger operators perturbed by potentials homogeneous with
respect to a certain set ${\mathfrak U}$ of unitary operators might
possess a lot of interesting properties. Obviously, such properties
became even more meaningful if, in addition to (\ref{e1}), the set
${\mathfrak U}$ has further algebraic group properties. In
particular, if ${\mathfrak U}$ is the set of scaling
transformations, then the additional multiplicative property
$U_{t_1}U_{t_2}=U_{t_2}U_{t_1}=U_{t_1t_2}$ of it elements enables
one to get simple solutions of many problems (like description of
nonnegative operator realizations, spectral properties, completeness
of the wave operators, explicit form of the scattering matrix) for
Schr\"{o}dinger operators with singular potentials
$\xi(t)$-invariant
 with respect to scaling transformations in ${\mathbb R^3}$
(Section~\ref{sec6}).

The abstract approach to the notion of symmetry developed in the
paper can be also useful for the study of supersingular
perturbations \cite{KP1}, for applications in the non-Archimedean
analysis (Example \ref{Examp}), and for the investigation of Weyl
families of boundary relations \cite{DHMS}.

In a very recent paper \cite{MT}, K. A. Makarov and E. Tsekanovskii
considered the so-called $\mu$-scale invariant operators, which can
be seen as a special case of $p(t)$-homogeneous operators in the
present paper. The main result of \cite{MT} is intimately related to
\cite[Lemma 4.5]{MS12}, see also Section 4 below.

Throughout the paper $\mathcal{D}(A)$, $\mathcal{R}(A)$, and
$\ker{A}$ denote the domain, the range, and the null-space of a
linear operator $A$, respectively, while
$A\upharpoonright{\mathcal{D}}$ stands for the restriction of $A$ to
the set $\mathcal{D}$.

\section{Preliminaries on operator realizations}
\label{sec2}

Following \cite{AL1, A4} an operator realization $A$ of (\ref{ne3})
in ${\mathfrak H}$ are defined by (\ref{lesia40}), (\ref{tat44}). To
clarify the meaning of $\mathbb{A}_0$ and $\psi_j^{\mathrm{ex}}$ in
(\ref{tat44}), observe that $\mathbb{A}_0$ stands for the
continuation of $A_0$ as a bounded linear operator acting from
${\mathfrak H}$ into ${\mathfrak H}_{-2}(A_0)$. Using the extended
resolvent $(\mathbb{A}_0+I)^{-1}$ this continuation can be
determined also by the formula
\begin{equation}\label{mmm5}
 \mathbb{A}_0f=[(\mathbb{A}_0+I)^{-1}]^{-1}f-f, \quad
 \forall{f}\in\mathfrak{H}.
\end{equation}
The linear functionals $<\psi_j^{\mathrm{ex}},\cdot>$ are extensions
of $<\psi_j,\cdot>$ onto ${\mathcal{D}}(A_{\mathrm{sym}}^*)$. Using
the well-known relation
\begin{equation}\label{k8}
\mathcal{D}(A_{\mathrm{sym}}^*)=\mathcal{D}(A_0)\dot{+}\mathcal{H},
\hspace{5mm} \mbox{where} \hspace{5mm} \mathcal{H}=\ker(A_{\mathrm{sym}}^*+I),
\end{equation}
one concludes that  $<\psi_j,\cdot>$ can be extended onto
$\mathcal{D}(A_{\mathrm{sym}}^*)$ by fixing their values on
$\mathcal{H}$. It follows from (\ref{ada5}) and (\ref{e7}) that the
vectors
\begin{equation}\label{kk41}
h_j=(\mathbb{A}_0+I)^{-1}\psi_j, \quad j=1,\ldots,n,
\end{equation}
form a basis of the defect subspace
$\mathcal{H}=\ker(A_{\mathrm{sym}}^*+I)$ of $A_{\mathrm{sym}}$.
Hence, the functionals $<\psi_j^{\mathrm{ex}},\cdot>$ are
well-defined by the formula
\begin{equation}\label{k23}
 <\psi_j^{\mathrm{ex}}, f>=<\psi_j, u>+\sum_{p=1}^n\alpha_p{r_{jp}}
\end{equation}
for all elements
$f=u+\sum_{p=1}^n\alpha_ph_p\in\mathcal{D}(A_{\mathrm{sym}}^*)$
$(u\in{\mathcal{D}}(A_0)$, $\alpha_p\in\mathbb{C})$ if the entries
$r_{jp}=<\psi_j, (\mathbb{A}_0+I)^{-1}\psi_p>=<\psi_j, h_p>$ of the
matrix $\mathbf{R}=(r_{jp})_{j,p=1}^n$ are known.

If all $\psi_j\in\mathfrak{H}_{-1}(A_0)$, then $r_{jp}$ are well
defined and $\mathbf{R}$ is a Hermitian matrix \cite{AL1}.
Otherwise, the matrix $\mathbf{R}$ is not uniquely determined. In
what follows, it is assumed that $\mathbf{R}$ is already chosen as a
Hermitian matrix. The problem of an appropriate choice of
$\mathbf{R}$ will be discussed in Section~\ref{sec3}.

In order to describe an operator realization $A$ of (\ref{ne3}) in
terms of parameters $b_{ij}$ of the singular perturbation $V$, the
method of boundary triplets (see \cite{DM, Gor} and the references
therein) is now incorporated.

\begin{definition}{\cite{Gor}}\label{d22}
A triplet $(N, \Gamma_0, \Gamma_1)$, where $N$ is an auxiliary
Hilbert space and $\Gamma_0$, $\Gamma_1$ are linear mappings of
$\mathcal{D}(A_{\mathrm{sym}}^*)$ into $N$, is called a boundary
triplet of $A_{\mathrm{sym}}^*$ if  $(A_{\mathrm{sym}}^*f, g)-(f,
A_{\mathrm{sym}}^*g)=(\Gamma_1f, \Gamma_0g)_{N}-(\Gamma_0f,
\Gamma_1g)_{N}$ for all $f, g\in\mathcal{D}(A_{\mathrm{sym}}^*)$ and
the mapping $(\Gamma_0, \Gamma_1) : \mathcal{D}(A_{\mathrm{sym}}^*)
\to N\oplus{N}$ is surjective.
\end{definition}

The next two results (Lemma \ref{l23} and Theorem \ref{ttt12}) are
known (see e.g. \cite{AKD, DHS}). For the convenience of the reader
some principal steps of their proofs are repeated.

\begin{lemma}\label{l23}
The triplet $({\mathbb C}^n, \Gamma_0, \Gamma_1)$, where the linear
operators $\Gamma_i:\mathcal{D}(A_{\mathrm{sym}}^*)\to{\mathbb C}^n$
are defined by the formulas
\begin{equation}\label{k9}
\Gamma_0f=\left(\begin{array}{c}
 <\psi_1^{\mathrm{ex}}, f> \\
 \vdots \\
 <\psi_n^{\mathrm{ex}}, f>
 \end{array}\right),
 \qquad \Gamma_1f=-\left(\begin{array}{c}
  \alpha_1 \\
  \vdots \\
 \alpha_n
  \end{array}\right),
\end{equation}
where $f=u+\sum_{j=1}\alpha_jh_j\in\mathcal{D}(A_{\mathrm{sym}}^*)$
\ $(u\in\mathcal{D}(A_0),\, \alpha_j\in\mathbb{C})$ and
$<\psi_j^{\mathrm{ex}}, f>$ is defined by (\ref{k23}), forms a
boundary triplet for $A_{\mathrm{sym}}^*$.
\end{lemma}

\begin{proof}
Using (\ref{ada5}), (\ref{k8}), and (\ref{kk41}) it is easy to
verify that the mappings
\begin{equation}\label{lesia99}
\widehat{\Gamma}_0f=\left(\begin{array}{c}
 \alpha_1 \\
 \vdots \\
\alpha_n
\end{array}\right), \qquad
\widehat{\Gamma}_1f=\left(\begin{array}{c}
 <\psi_1, u> \\
 \vdots \\
<\psi_n, u>
\end{array}\right), \quad f=u+\sum_{j=1}\alpha_jh_j
\end{equation}
satisfy the conditions of Definition \ref{d22}. Thus $({\mathbb
C}^n, \widehat{\Gamma}_0, \widehat{\Gamma}_1)$ is a boundary triplet
for $A_{\mathrm{sym}}^*$. It follows from (\ref{k23}), (\ref{k9}),
and (\ref{lesia99}) that
\begin{equation}\label{lesia101}
\Gamma_0f=\widehat{\Gamma}_1f+{\mathbf R}\widehat{\Gamma}_0f, \quad
\Gamma_1f=-\widehat{\Gamma}_0f, \qquad
f\in\mathcal{D}(A_{\mathrm{sym}}^*).
\end{equation}
These relations between ${\Gamma}_i$ and $\widehat{\Gamma}_i$ and
the fact that $({\mathbb C}^n, \widehat{\Gamma}_0,
\widehat{\Gamma}_1)$ is a boundary triplet for $A_{\mathrm{sym}}^*$
imply that $({\mathbb C}^n, \Gamma_0, \Gamma_1)$ also is a boundary
triplet for $A_{\mathrm{sym}}^*$.
\end{proof}

\begin{theorem}\label{ttt12}
The operator realization $A$ of (\ref{ne3}) is an intermediate
extension of $A_{\mathrm{sym}}$ which coincides with the operator
\begin{equation}\label{k41}
 A_{\mathbf{B}}=A_{\mathrm{sym}}^*\upharpoonright{\mathcal{D}(A_{\mathbf{B}})},
 \quad
 {\mathcal{D}(A_{\mathbf{B}})}=\{\, f\in\mathcal{D}(A_{\mathrm{sym}}^*) :\,
 {\mathbf B}\Gamma_0f=\Gamma_1f \,\},
\end{equation}
where $\Gamma_i$ are defined by (\ref{k9}) and
$\mathbf{B}=(b_{ij})_{i,j=1}^n$ is the coefficient matrix of the
singular perturbation
$V=\sum_{i,j=1}^{n}{b}_{ij}<\psi_j,\cdot>\psi_i$ in (\ref{ne3}).

If $V$ is symmetric, i.e., $<Vu, v>=<u, Vv>$
$(u,v\in\mathfrak{H}_2(A_0))$, then the corresponding operator
realization $A_{\mathbf{B}}$ becomes self-adjoint.
\end{theorem}
\begin{proof}
It follows from (\ref{mmm5}) that $\mathbb{A}_0h_j=\psi_j-h_j$ for
all $h_j$ defined by (\ref{kk41}). Rewriting
$f\in\mathcal{D}(A_{\mathrm{sym}}^*)$ in the form
$f=u+\sum_{i=1}\alpha_ih_i$, where $u\in{\mathcal D}(A_0)$,
$h_i\in{\mathcal H}$, $\alpha_i\in\mathbb{C}$, and using
(\ref{tat44}) and (\ref{k9}) leads to
\[
\begin{split}
 \mathbb{A}_\mathbf{R}f
 &=A_0u-\sum_{i=1}^{n}\alpha_ih_i
   +\sum_{i,j=1}^n{b}_{ij}<\psi_j^{\mathrm{ex}},f>\psi_i+
 \sum_{i=1}^n\alpha_i\psi_i  \\
 &=A_{\mathrm{sym}}^*f+(\psi_1,\ldots,\psi_n)[{\mathbf B}\Gamma_0f-\Gamma_1f].
\end{split}
\]
This equality and (\ref{lesia40}) show that $f\in\mathcal{D}(A)$ if
and only if  ${\mathbf B}\Gamma_0f-\Gamma_1f=0$. Therefore, the
operator realization $A$ of (\ref{ne3}) is an intermediate extension
of $A_{\mathrm{sym}}$ and $A$ coincides with the operator
$A_{\mathbf{B}}$ defined by (\ref{k41}).

To complete the proof it suffices to finally observe that $V$ is
symmetric if and only if the corresponding matrix of coefficients
$\mathbf{B}=(b_{ij})_{i,j=1}^n$ is Hermitian. In this case
(\ref{k41}) immediately implies the self-adjointness of
$A_{\mathbf{B}}$.
\end{proof}

\begin{remark}
Another approach, also involving the use of boundary triplets, to
determine self-adjoint operator realizations of finite rank singular
perturbations of the form $A_0+G\alpha{G^*}$, where $G$ is an
injective linear mapping from $\mathbb{C}^n$ to
${\mathfrak{H}}_{-k}(A_0)$ was presented in \cite[Section 4]{DHS}.
\end{remark}

\section{Admissible matrices and admissible operators}
\label{sec3}

There are certain natural requirements for the determination of the
entries $r_{jp}$ of the matrix $\mathbf{R}$ in (\ref{k23}). Indeed,
if the linear span $\mathcal{X}$ of $\{\psi_j\}_{j=1}^n$ has a
nonzero intersection with $\mathfrak{H}_{-1}(A_0)$, then for any
$\psi\in\mathcal{X}\cap\mathfrak{H}_{-1}(A_0)$, the corresponding
element $h=(\mathbb{A}_0+I)^{-1}\psi$ belongs to
$\mathfrak{H}_{1}(A_0)$ and, hence, the functional $<\psi, \cdot>$
defined by (\ref{ada5}) admits the following extension by continuity
onto $\mathfrak{H}_{1}(A_0)$:
$$
<\psi,f>=((A_0+I)^{1/2}f,(A_0+I)^{1/2}h), \hspace{5mm}
\forall{f}\in{\mathfrak
 H}_{1}(A_0).
$$
To preserve such natural extensions of $<\psi, \cdot>$ onto
${\mathcal{D}}(A_{\mathrm{sym}}^*)\cap{\mathfrak{H}}_{1}(A_0)$ in
the definition (\ref{k23}), the concept of admissible matrices
$\mathbf{R}$ as introduced in \cite{A4} is used.

\begin{definition}\label{d11}
A Hermitian matrix $\mathbf{R}=(r_{jp})_{j,p=1}^n$ is called
admissible for the regularization $\mathbb{A}_\mathbf{R}$ of
(\ref{ne3}) if its entries $r_{jp}$ are chosen in such a way that if
a singular element $\psi=c_1\psi_1+\cdots +c_n\psi_n$ belongs to
$\mathfrak{H}_{-1}(A_0)$, then for all
$f\in\mathcal{D}(A_{\mathrm{sym}}^*)\cap{\mathfrak H}_{1}(A_0)$
\begin{equation}\label{les21}
 <\psi^{\mathrm{ex}},f>=((A_0+I)^{1/2}f,(A_0+I)^{1/2}h)=\sum_{j=1}^n{c_j}<\psi_j^{\mathrm{ex}}, f>,
\end{equation}
where $<\psi_j^{\mathrm{ex}},f>$ are defined by (\ref{k23}) and
$h=(\mathbb{A}_0+I)^{-1}\psi$.
\end{definition}

It is convenient to describe the set of admissible matrices in terms
of a certain associated operators. It follows from relations
(\ref{lesia101}) that the choice of a matrix ${\mathbf R}$ in
(\ref{k23}) \textit{is equivalent} to the choice of an operator
$\widetilde{A}$ defined by
\begin{equation}\label{lesia200}
 \widetilde{A}:=A_{\mathrm{sym}}^*\upharpoonright{\mathcal{D}(\widetilde{A})},
 \quad \mathcal{D}(\widetilde{A})
  =\ker\Gamma_0=\{\, f\in{\mathcal D}(A_{\mathrm{sym}}^*)  :\,
  -{\mathbf R}\widehat{\Gamma}_0f=\widehat{\Gamma}_1f \,\}.
\end{equation}

\begin{definition}\label{d1}
An operator $\widetilde{A}$ is called admissible for the
regularization  of (\ref{ne3}) if $\widetilde{A}$ is defined by
(\ref{lesia200}) with an admissible matrix ${\mathbf R}$.
\end{definition}

Since ${\mathbf R}$ is Hermitian, Definition~\ref{d1} and the
general theory of boundary triplets \cite{DM} imply that an
admissible operator $\widetilde{A}$ is a self-adjoint extension of
$A_{\mathrm{sym}}$. In general, $\widetilde{A}$ need not be
nonnegative. It is nonnegative if and only if
\begin{equation}\label{sas44}
 (A_F+I)^{-1}\leq(\widetilde{A}+I)^{-1}\leq(A_N+I)^{-1},
\end{equation}
where  $A_F$ is the Friedrichs extension and  $A_N$ is the Krein-von
Neumann extension of $A_{\mathrm{sym}}$ (see e.g., \cite{HMS} and
the references therein).

The next lemma gives some useful facts concerning the (unperturbed)
nonnegative self-adjoint operator $A_0$ and its relation to the
Friedrichs extension $A_F$ of $A_{\mathrm{sym}}$. They can be
considered to be well known from the extension theory of nonnegative
operators, therefore details for the present formulations with their
proofs are left to the reader; see e.g. \cite{AN, DM95, HMS, HSSW,
Krein, KK}.
\begin{lemma}\label{newle}
Let $C=(A_0+I)^{-1}-(A_F+I)^{-1}$ and let $S_0=A_0\cap A_F$.
Moreover, denote $\cH=\ker(A_{\mathrm{sym}}^*+I)$ and
$\cH'=\ker(S_0^*+I)$. Then:
\begin{enumerate}
\def\labelenumi{\rm (\roman{enumi})}
\item $\overline{\cR(C)}=\cH'$;
\item $\ker C=\cR(S_0+I)= \cR(A_{\mathrm{sym}}+I)\oplus \cH''$,
where $\cH''=\cH\ominus \cH'$;
\item $\cR(C^{1/2})=\cD(A_0^{1/2})\cap \cH={\cH'}$; %\cD(A_0^{1/2})\cap \cH'$;
\item $\cD(A_0^{1/2})=\cD(A_F^{1/2}) \dot{+} \cR(C^{1/2})$.
\end{enumerate}
\end{lemma}

Using the spaces introduced in \eqref{ee1} and (iii) in Lemma
\ref{newle} one can rewrite the decomposition in part (iv) of
Lemma~\ref{newle} as follows:
\begin{equation}\label{kk2}
\mathfrak{H}_{1}(A_0)=\mathcal{D}\oplus_1{\mathcal{H}}', \quad
 {\mathcal{H}}'={\mathcal{H}}\cap\mathfrak{H}_{1}(A_0)
 =(\mathbb{A}_0+I)^{-1}[\mathcal{X}\cap\mathfrak{H}_{-1}(A_0)],
\end{equation}
where $\mathcal{D}\,(=\cD(A_F^{1/2}))$ stands for the completion of
$\mathcal{D}(A_{\mathrm{sym}})$ in $\mathfrak{H}_{1}(A_0)$,
$\oplus_1$ denotes the orthogonal sum in $\mathfrak{H}_{1}(A_0)$,
and $\mathcal{X}$ is the linear span of $\{\psi_j\}_{j=1}^n$.

The set of all admissible operators can now be characterized in
'coordinate free' manner as follows.

\begin{theorem}\label{td1}
A self-adjoint extension $\widetilde{A}$ of $A_{\mathrm{sym}}$ is an
admissible operator for the regularization of (\ref{ne3}) if and
only if  $\widetilde{A}$ is transversal to $A_0$ (i.e.,
$\mathcal{D}(A_0)+\mathcal{D}(\widetilde{A})=\mathcal{D}(A_{\mathrm{sym}}^*))$
and
\begin{equation}\label{kak12}
 \mathcal{D}(\widetilde{A})\cap{\mathfrak
 H}_{1}(A_0)\subset{\mathcal{D}(A_F)},
\end{equation}
where $A_F$ is the Friedrichs extension of $A_{\mathrm{sym}}$.
\end{theorem}
\begin{proof}
Assume that the self-adjoint extension $\widetilde{A}$ of
$A_{\mathrm{sym}}$ is transversal to $A_0$ and it satisfies the
condition (\ref{kak12}). In view of (\ref{lesia99}),
$\mathcal{D}({A}_0)=\ker\widehat{\Gamma}_0$. Therefore, the
transversality of $\widetilde{A}$ and $A_0$ is equivalent to the
representation of $\mathcal{D}(\widetilde{A})$ in the form
(\ref{lesia200}) with an $n\times n$ Hermitian matrix ${\mathbf R}$
(here $A_{\mathrm{sym}}$ has finite defect numbers $(n,n)$), cf.
\cite[Proposition~1.4]{DM95}.

Since
\begin{equation} \label{DAF}
 \mathcal{D}(A_F)=\mathcal{D}\cap\mathcal{D}(A_{\mathrm{sym}}^*),
\end{equation}
the decomposition (\ref{kk2}) shows that the condition (\ref{kak12})
is equivalent to the relation
\begin{equation}\label{kk12}
 ((A_0+I)^{1/2}\tilde{f},(A_0+I)^{1/2}h)=0, \quad
  \forall{\tilde{f}\in\mathcal{D}(\widetilde{A})\cap{\mathfrak H}_{1}(A_0)}, \quad
  \forall{h}\in{\mathcal{H}}'.
\end{equation}
Now it is shown that ${\mathbf R}$ is an admissible matrix in the
sense of Definition~\ref{d11} by verifying (\ref{les21}) for all
$\psi\in\mathcal{X}\cap\mathfrak{H}_{-1}(A_0)$. Observe, that the
mapping $\Gamma_0$ defined in Lemma~\ref{l23}, see also
(\ref{lesia101}), determines the extended functionals
$<\psi_j^{\mathrm{ex}}, f>$ in (\ref{k23}).

The transversality of $\widetilde{A}$ and $A_0$ yields the following
decomposition for the elements
$f\in\mathcal{D}(A_{\mathrm{sym}}^*)$:
\begin{equation}\label{kk14}
 f=\tilde{f}+u,
\end{equation}
where $\tilde{f}\in\mathcal{D}(\widetilde{A})$ and
$u\in\mathcal{D}(A_0)$ are uniquely determined modulo
$\mathcal{D}(A_{\mathrm{sym}})$. If
$\psi=\sum_{j=1}^nc_j\psi_j\in\mathfrak{H}_{-1}(A_0)$, then by
(\ref{kk2}) $h=(\mathbb{A}_0+I)^{-1}\psi\in{\mathcal{H}}'$. Now with
$f\in\mathcal{D}(A_{\mathrm{sym}}^*)\cap{\mathfrak H}_{1}(A_0)$
decomposed as in (\ref{kk14}) one obtains:
\begin{eqnarray}
\label{eqfunc}
 <\psi^{\mathrm{ex}},f>\,\,=\,\sum_{j=1}^nc_j<\psi_j^{\mathrm{ex}},f>\,
 =\mathbf{c}{\Gamma}_0{f}\stackrel{(\ref{kk14})}{=}
  \mathbf{c}{\Gamma}_0(\tilde{f}+u) & &  \\ \nonumber
 \stackrel{(\ref{lesia101})}{=}
 \mathbf{c}(\widehat{\Gamma}_1+{\mathbf R}\widehat{\Gamma}_0)u
 =\mathbf{c}\widehat{\Gamma}_1 u
 \stackrel{(\ref{lesia99})}{=}<\psi,u>
 \stackrel{(\ref{ada5})}{=}((A_0+I)u,h)
\end{eqnarray}
where $\mathbf{c}:=(c_1,\ldots,c_n)$. On the other hand, it follows
from (\ref{kk12}) that
$$
 ((A_0+I)^{1/2}f,(A_0+I)^{1/2}h)=
 ((A_0+I)^{1/2}(\tilde{f}+u),(A_0+I)^{1/2}h)
 =((A_0+I)u,h),
$$
which combined with (\ref{eqfunc}) proves (\ref{les21}). Thus,
${\mathbf R}$ is an admissible matrix and $\widetilde{A}$ is an
admissible operator.

Conversely, assume that $\widetilde{A}$ is an admissible operator.
Then the relation (\ref{lesia200}) ensures the transversality of
$\widetilde{A}$ and $A_0$ and ${\mathbf R}$ determines the extended
functionals $<\psi_j^{\mathrm{ex}}, \cdot>$  via (\ref{k23}).
Reasoning as in (\ref{eqfunc}) it is seen that (\ref{les21}) implies
$$
 0=((A_0+I)^{1/2}f,(A_0+I)^{1/2}h)-<\psi^{\mathrm{ex}},f>\,\,
 =((A_0+I)^{1/2}\tilde{f},(A_0+I)^{1/2}h)
$$
for all ${f\in\mathcal{D}(A_{\mathrm{sym}}^*)\cap{\mathfrak
H}_{1}(A_0)}$ and ${h}\in{\mathcal{H}}'$. Thus, the relation
(\ref{kk12}) and, equivalently, the relation (\ref{kak12}) is
satisfied. Theorem~\ref{td1} is proved.
\end{proof}

For some further study of admissible operators the following lemma
is needed.

\begin{lemma}\label{lil1}
Let $\widetilde{\mathcal H}$ be a subspace of
${\mathcal{H}}=\ker(A_{\mathrm{sym}}^*+I)$. Then the symmetric
operator
\begin{equation}\label{les36}
 S=A_F\upharpoonright_{\mathcal{D}(S)}, \quad
 \mathcal{D}(S)
 =(A_F+I)^{-1}[\mathcal{R}(A_{\mathrm{sym}}+I)\oplus{\widetilde{\mathcal{H}}}]
\end{equation}
 satisfies the relations
\begin{equation}\label{sas57}
 \mathcal{D}(S)\cap\mathcal{D}(A_0)=\mathcal{D}(A_{\mathrm{sym}}) \quad
 \mbox{and} \quad \mathcal{D}(S)+\mathcal{D}(A_0)
 =\mathcal{D}(A_F)\dot{+}{\mathcal{H}}'
\end{equation}
if and only if
\begin{equation}\label{dur57}
\dim\widetilde{\mathcal{H}}=\dim{\mathcal{H}}' \quad \mbox{and}
\quad \widetilde{\mathcal{H}}\cap{\mathcal{H}}''=\{0\},
\end{equation}
where $\cH'=\cH\cap\mathfrak{H}_{1}(A_0)$ and
$\cH''=\cH\ominus\cH'$. In this case, the domain of $S$ admits the
description
\begin{equation}\label{sas90}
 \mathcal{D}(S)=\mathcal{D}(A_{\mathrm{sym}})\,\dot{+}\,\{\, h'+u :\,
  {h'}\in{\mathcal{H}}',\quad u=u(h') \,\},
\end{equation}
where $u=u(h')\in\mathcal{D}(A_0)$ is (uniquely) determined by
$h'\in\cH'$ and satisfies the relation
\begin{equation}\label{sas70}
((A_0+I)u, \wt{h}^{\perp})=<\psi, u> \, =0, \quad
\forall\wt{h}^{\perp}\in{\mathcal{H}}\ominus\widetilde{{\mathcal{H}}},
\quad \psi=(\mathbb{A}_0+I)\wt{h}^{\perp}.
\end{equation}
\end{lemma}

\begin{proof} Denote
$S_0=A_F\cap A_0$. By Lemma~\ref{newle}
\begin{equation}\label{tut}
 \cD(S_0)=(A_0+I)^{-1}[\mathcal{R}(A_{\mathrm{sym}}+I)\oplus\mathcal{H}'']
 =(A_F+I)^{-1}[\mathcal{R}(A_{\mathrm{sym}}+I)\oplus\mathcal{H}''],
\end{equation}
where ${\mathcal{H}}''={\mathcal{H}}\ominus{\mathcal{H}}'$.
Comparing (\ref{les36}) and (\ref{tut}), one concludes that
$$
 \mathcal{D}(S)\cap\mathcal{D}(A_0)
 =\mathcal{D}(S)\cap\mathcal{D}(S_0)
 =(A_F+I)^{-1}[\mathcal{R}(A_{\mathrm{sym}}+I)
    \oplus(\widetilde{\mathcal{H}}\cap\mathcal{H}'')].
$$
Thus,
$$
 \mathcal{D}(S)\cap\mathcal{D}(A_0)=\mathcal{D}(A_{\mathrm{sym}})
\iff\widetilde{\mathcal{H}}\cap\mathcal{H}''=\{0\}.
$$
The relations (\ref{les36}) and (\ref{tut}) also show that
\begin{equation}\label{tut1}
 \mathcal{D}(S)+\mathcal{D}(A_0)
 =(A_F+I)^{-1}[\mathcal{R}(A_{\mathrm{sym}}+I)
   \oplus(\widetilde{\mathcal{H}}\dot{+}\mathcal{H}'')]
   +(A_0+I)^{-1}{\mathcal{H}}'.
\end{equation}
Here $(A_0+I)^{-1}{\mathcal{H}}'$  can be represented as
\begin{equation}\label{tut4}
 (A_0+I)^{-1}{\mathcal{H}}'=\{\,(A_F+I)^{-1}h'+Ch' :\,
     h' \in {\mathcal{H}}' \,\},
\end{equation}
where $C=(A_0+I)^{-1}-(A_F+I)^{-1}$. It follows from
Lemma~\ref{newle} that
\begin{equation}\label{tut6}
  \cR(C)=\mathcal{H}', \quad  \ker{C}=\ran(A_{\mathrm{sym}}+I)\oplus \cH''.
\end{equation}
Relations (\ref{tut1}), (\ref{tut4}), and (\ref{tut6}) show that the
second identity in (\ref{sas57}) holds  if and only if
$\widetilde{\mathcal{H}}\dot{+}\mathcal{H}''=\mathcal{H}$.
Obviously, this representation is possible only in the case where
$\dim\widetilde{\mathcal{H}}=\dim{\mathcal{H}}'$.

The definition (\ref{les36}) shows that
$\mathcal{D}(S)=\mathcal{D}(A_{\mathrm{sym}})\dot{+}(A_F+I)^{-1}\widetilde{\mathcal{H}},$
where
$$
 (A_F+I)^{-1}\widetilde{\mathcal{H}}=\{\,(A_0+I)^{-1}\wt{h}-C\wt{h}
 :\,\, {\wt h}\in\widetilde{\mathcal{H}} \,\}.
$$
Since $\widetilde{\mathcal{H}}$ satisfies (\ref{dur57}), it follows
from (\ref{tut6}) that $C\widetilde{\mathcal{H}}={\mathcal{H}}'$.
Now, setting $u=(A_0+I)^{-1}\wt{h}$ and $h'=-C\wt{h}$, one obtains
(\ref{sas90}) and (\ref{sas70}). Note that the preimage
$\wt{h}=C^{-1}h'\in\wt\cH$, and therefore also $u$, is uniquely
determined by $h'\in\cH'$,
\end{proof}

The next theorem gives a description of all admissible operators.

\begin{theorem}\label{did121}
Let $\wt A$ be a self-adjoint extension of $A_{\mathrm{sym}}$ and
let the symmetric operator $S=\wt A\cap{A_F}$ be represented as in
\eqref{les36} with some subspace $\wt\cH$ of $\cH$. Then the
following statements are equivalent:
\begin{enumerate}
\def\labelenumi{\rm (\roman{enumi})}
\item $\wt A$ is an admissible operator for the regularization of
\eqref{ne3};

\item $\wt A$ is a self-adjoint extension of
$S$ transversal to the Friedrichs extension $S_F$ of $S$ and the
subspace $\widetilde{\mathcal{H}}$ satisfies the conditions in
\eqref{dur57}.
\end{enumerate}
\end{theorem}

\begin{proof}
Let $\wt A$ be an admissible operator. Since $\wt A$ and $A_0$ are
transversal, one has
\begin{equation}\label{tut3}
\mathcal{D}(\wt
A)\cap\mathcal{D}(A_0)=\mathcal{D}(A_{\mathrm{sym}}), \quad
\mathcal{D}(\wt
A)+\mathcal{D}(A_0)=\mathcal{D}(A_F)\dot{+}{\mathcal{H}}=\mathcal{D}(A^*_{\mathrm{sym}}).
\end{equation}
The condition \eqref{kak12} is equivalent to
$$
 \cD(\wt A) \cap \sH_1(A_0)=\cD(\wt A) \cap \cD(A_F)=\cD(\wt A\cap
 A_F).
$$
Thus, intersecting all parts of (\ref{tut3}) with ${\mathfrak
H}_{1}(A_0)$ one concludes that the relations (\ref{sas57}) are true
for $S=\wt A\cap A_F$. By Lemma~\ref{lil1}, the subspace
$\widetilde{\mathcal H}$ satisfies (\ref{dur57}). Furthermore, since
the Friedrichs extension $S_F$ of $S$ coincides with $A_F$, one gets
$\mathcal{D}(S_F)\cap\mathcal{D}(\widetilde{A})
=\mathcal{D}(A_F)\cap\mathcal{D}(\widetilde{A})=\mathcal{D}(S)$.
This implies the transversality of $S_F$ and $\wt A$. The
implication (i) $\Rightarrow$ (ii) is proved.

Now, assume that (ii) is satisfied. Since
$S\supset{A_{\mathrm{sym}}}$, the operator $\widetilde{A}$ is a
self-adjoint extension of $A_{\mathrm{sym}}$. It follows from
(\ref{les36}) that $\ker(S^*+I)={\mathcal
H}\ominus{\widetilde{\mathcal H}}$ and hence,
$\mathcal{D}(S^*)=\mathcal{D}(S_F)+\ker(S^*+I)=\mathcal{D}(A_F)\dot{+}({\mathcal
H}\ominus{\widetilde{\mathcal H}})$. On the other hand, the
transversality of $S_F$ and $\wt A$ gives
$\mathcal{D}(S^*)=\mathcal{D}(A_F)+\mathcal{D}(\widetilde{A})$.
Therefore,
$\mathcal{D}(A_F)+\mathcal{D}(\widetilde{A})=\mathcal{D}(A_F)\dot{+}({\mathcal
H}\ominus{\widetilde{\mathcal H}})$. This equality and the second
relation in (\ref{sas57}) yield
\begin{equation}
\label{eq2}
\begin{split}
 \mathcal{D}(A_0)+\mathcal{D}(\widetilde{A})
 &=\mathcal{D}(S)+\mathcal{D}(A_0)+\mathcal{D}(\widetilde{A}) \\
 &=(\mathcal{D}(A_F)\dot{+}{\mathcal{H}}')+\mathcal{D}(\widetilde{A})=\mathcal{D}(A_F)\dot{+}{\mathcal{H}'}\dot{+}({\mathcal H}
   \ominus{\widetilde{\mathcal H}}).
\end{split}
\end{equation}
The conditions (\ref{dur57}) imply that
${\mathcal{H}'}\dot{+}({\mathcal H}\ominus{\widetilde{\mathcal
H}})={\mathcal{H}}$. Hence, (\ref{eq2}) shows that
$\mathcal{D}(A_0)+\mathcal{D}(\widetilde{A})
=\mathcal{D}(A_F)\dot{+}{\mathcal{H}}=\mathcal{D}(A^*_{\mathrm{sym}})$,
i.e., $\widetilde{A}$ and $A_0$ are transversal. Furthermore, by
Lemma~\ref{newle}, see also (\ref{DAF}),
$\mathcal{D}(A_F)\dot{+}{\mathcal{H}}'
 =\mathfrak{H}_{1}(A_0)\cap\mathcal{D}(A_{\mathrm{sym}}^*)$.
Now, employing the second relation in (\ref{sas57}) one obtains
$$
 \mathcal{D}(\widetilde{A})\cap\mathfrak{H}_{1}(A_0)
 =\mathcal{D}(\widetilde{A})\cap(\mathcal{D}(S)+\mathcal{D}(A_0))
 =\mathcal{D}(S){+}\mathcal{D}(A_{\mathrm{sym}})
 =\mathcal{D}(S)\subset{\mathcal{D}(A_{F})}.
$$
According to Theorem \ref{td1} this means that $\widetilde{A}$ is an
admissible operator for the regularization of (\ref{ne3}). Thus, the
implication (ii) $\Rightarrow$ (i) is proved.
\end{proof}

It follows from Theorem \ref{did121} that there is at least one
admissible operator for the regularization of (\ref{ne3}).

\begin{cor}\label{p121}
If all the elements $\psi_j$ in (\ref{ne3}) belong to
$\mathfrak{H}_{-1}(A_0)$, then there exists a unique admissible
operator for the regularization of (\ref{ne3}) and it coincides with
the Friedrichs extension $A_F$ of $A_{\mathrm{sym}}$.
\end{cor}
\begin{proof} Assume that $\psi_j\in\mathfrak{H}_{-1}(A_0)$ for all $j=1,\dots,n$.
Then $\mathcal{D}(A_{\mathrm{sym}}^*)\subset\mathfrak{H}_{1}(A_0)$
and $\cH'=\cH$. Let $\widetilde{A}$ be an admissible operator for
the regularization of (\ref{ne3}) and let $S=\wt A\cap A_F$. By
Theorem~\ref{did121} the corresponding subspace $\wt\cH$ satisfies
(\ref{dur57}) in Lemma~\ref{lil1}, so that $\wt\cH=\cH$. Now
(\ref{les36}) gives $S=A_F$ and since $S=\wt A\cap A_F$, one
concludes that $\widetilde{A}=A_F$. This completes the proof.
\end{proof}

\begin{cor}\label{p1211}
If all the elements $\psi_j$ in (\ref{ne3}) are
$\mathfrak{H}_{-1}(A_0)$-independent (i.e.
$\mathcal{X}\cap\mathfrak{H}_{-1}(A_0)=\{0\}$), then every
self-adjoint extension $\widetilde{A}$ of $A_{\mathrm{sym}}$
transversal to $A_0$ is admissible for the regularization of
(\ref{ne3}). The Friedrichs extension of $A_{\mathrm{sym}}$
coincides with $A_0$.
\end{cor}
\begin{proof} The condition of $\mathfrak{H}_{-1}(A_0)$-independency
means that ${\mathcal{H}}'=\{0\}$. In this case, only the zero
subspace $\widetilde{\mathcal H}=\{0\}$ can satisfy (\ref{dur57}).
The corresponding operator $S$ coincides with $A_{\mathrm{sym}}$.
Moreover, since ${\mathcal{H}}'=\{0\}$, Lemma~\ref{newle} shows that
$S_F=A_F=A_0$. Thus, by Theorem \ref{did121}, $\wt A$ is admissible
if and only if it is transversal to $A_0$.
\end{proof}

The properties of admissible operators are closely related to the
transversality of the Friedrichs and the Krein-von Neumann
extensions of $A_{\mathrm{sym}}$.

\begin{theorem}\label{t303}
There exists a nonnegative admissible operator $\wt A$ for the
regularization of (\ref{ne3}) if and only if the Friedrichs
extension $A_F$ and the Krein-von Neumann extension $A_N$ of
$A_{\mathrm{sym}}$ are transversal.
\end{theorem}

\begin{proof} Let $\wt A$ be a nonnegative admissible operator.
Then $\widetilde{A}$ is a nonnegative extension of
$A_{\mathrm{sym}}$ and therefore $(\widetilde{A}+I)^{-1}$ satisfies
the inequalities \eqref{sas44}. Recall that transversality of
self-adjoint extensions $\widetilde{A}_1$ and $\wt A_2$ of
$A_{\mathrm{sym}}$ is equivalent to
\begin{equation}\label{sas404}
  [(\wt A_1+I)^{-1}-(\wt A_2+I)^{-1}]{\mathcal H}={\mathcal H}
\end{equation}
(see e.g. \cite{DM}). Hence, if $A_F$ and $A_N$ are not transversal
then $(A_F+I)^{-1}h=(A_N+I)^{-1}h$ for some nonzero $h\in{\mathcal
H}$. Then nonnegativity of $\widetilde{A}$ and $A_0$ yields
$(\widetilde{A}+I)^{-1}h=(A_0+I)^{-1}h$ due to (\ref{sas44}) (with
similar inequalities for $A_0$), so that
$$
[(\widetilde{A}+I)^{-1}-(A_0+I)^{-1}]{\mathcal H}\subset{\mathcal
H}\ominus<h>
$$
and by (\ref{sas404}) $\wt A$ and $A_0$ cannot be transversal. This
is a contradiction to the admissibility of $\wt A$. Thus $A_F$ and
$A_N$ are transversal.

To prove the converse statement assume that $A_F$ and $A_N$ are
transversal. Let $\wt\cH$ be a subspace of $\cH$, which satisfies
(\ref{dur57}) and let the symmetric operator $S$ be defined by
(\ref{les36}) in Lemma~\ref{lil1}. Moreover, let $\widetilde{A}$ be
the Krein-von Neumann extension of $S$. Clearly, $\wt A$ is a
nonnegative self-adjoint extension of $A_{\mathrm{sym}}$. It remains
to prove that the operator $\widetilde{A}$ is admissible for the
regularization of (\ref{ne3}). To see this, observe that the
Friedrichs extension of $S$ coincides with $A_F$. Then it follows
from \cite[Proposition~7.2]{AHSS} that the Friedrichs extension
$S_F=A_F$ and the Krein-von Neumann extension $\widetilde{A}$ of $S$
are transversal with respect to $S$. Therefore, by Theorem
\ref{did121}, $\widetilde{A}$ is an admissible operator.
\end{proof}

Observe that $S$ in Theorem~\ref{t303} is a restriction of the
Friedrichs extension $A_F$ of $A_{\mathrm{sym}}$. Since the
admissible operator $\widetilde{A}$ constructed in
Theorem~\ref{t303} is the Krein-von Neumann extension of $S$ it is a
consequence of \cite[Theorem 6.4]{AHSS} that $\widetilde{A}$ is an
extremal extension of $A_{\mathrm{sym}}$ in the sense of the
following definition
\begin{definition}{\cite{Ar, AHSS}}\label{d5}
A self-adjoint extension $\widetilde{A}$ of $A_{\mathrm{sym}}$ is
called extremal if it is nonnegative and satisfies the condition
$$
 \inf_{u\in{\mathcal{D}(A_{\mathrm{sym}})}}(\widetilde{A}(f-u), f-u)=0
 \quad \mathrm{for \ all} \ \ f\in\mathcal{D}(\widetilde{A}).
$$
\end{definition}

\begin{theorem}\label{did3}
Let the Friedrichs extension $A_F$ and the Krein-von Neumann
extension $A_N$ of $A_{\mathrm{sym}}$ be transversal, and let $S$ be
defined by (\ref{les36}) and (\ref{dur57}). Then among all
self-adjoint extensions of $S$ there exists a unique extremal
admissible operator $\widetilde{A}$ for the regularization of
(\ref{ne3}).
\end{theorem}
\begin{proof}
In view of Theorem \ref{t303}, it suffices to show that the
Krein-von Neumann extension $\widetilde{A}$ of $S$ is the only
extremal extension of $A_{\mathrm{sym}}$ which is admissible for the
regularization of (\ref{ne3}).

To prove this assume that $\widehat{A}$ is extremal and admissible.
Then by \cite[Theorem 6.4]{AHSS} $\widehat{A}$ as an extremal
extension of $A_{\mathrm{sym}}$ is the Krein-von Neumann extension
of the symmetric operator $\widehat{S}:=\widehat{A}\cap{A_F}$.
Moreover, by Theorem \ref{did121} the admissibility of $\widehat{A}$
means that $\widehat{S}$ is determined via (\ref{les36}) where the
corresponding subspace $\widehat{\mathcal H}$ satisfies
(\ref{dur57}).

Since $\widehat{A}$ is an extension of $S$, one has
$S\subseteq{\widehat{S}}$ or, equivalently, $\widetilde{\mathcal
H}\subseteq\widehat{\mathcal H}$, where the subspaces
$\widetilde{\mathcal H}$ and $\widehat{\mathcal H}$ correspond to
$S$ and ${\widehat{S}}$ in (\ref{les36}). Now the first equality in
(\ref{dur57}) forces that $\widetilde{\mathcal H}=\widehat{\mathcal
H}$ and hence $S=\widehat{S}$. Therefore,
$\widehat{A}=\widetilde{A}$ and this completes the proof.
\end{proof}

\begin{remark}\label{r1}
The selection of a self-adjoint operator $\widetilde{A}$ transversal
to the initial one $A_0$ (but without the admissibility condition
(\ref{kak12})) is also a key point of the approach used in
\cite{ARC} to the determination of self-adjoint realizations of a
formal expression $A_0+V$, where a singular perturbation $V$ is
assumed to be (in general) an unbounded self-adjoint operator $V:
\mathfrak{H}_2(A_0)\to\mathfrak{H}_{-2}(A_0)$ such that $\ker{V}$ is
dense in $\mathfrak{H}$.  In this case, the regularization of
$A_0+V$ takes the form $A_{\mathcal{P},
V}=\mathbb{A}_0+V{\mathcal{P}}$ and it is well defined on the domain
$ \mathcal{D}(A_{\mathcal{P},V})
 =\{\,f\in\mathcal{D}(A_{\mathrm{sym}}^*)  :\,
  {\mathcal P}f\in\mathcal{D}(V)\,\},$
where ${\mathcal P}$ is the skew projection onto
$\mathfrak{H}_{2}(A_0)$ in ${\mathcal D}(A_{\mathrm{sym}}^*)$ that
is uniquely determined by the choice of $\widetilde{A}$.
\end{remark}

\section{Singular perturbations with symmetries and uniqueness of admissible operators}
\label{sec4}

According to \eqref{k23} and \eqref{lesia200} the regularization
$\mathbb{A}_\mathbf{R}$ of (\ref{ne3}) depends on the choice of an
admissible operator $\widetilde{A}$. Apart from the case of form
bounded singular perturbations, admissible operators are not
determined uniquely, cf. Theorem \ref{did121}. However, in many
cases (see e.g., \cite{AL1, A4}), the uniqueness can be attained by
imposing extra assumptions of symmetry motivated by the specific
nature of the underlying physical problem. In this section, we study
this problem in an abstract framework.

\subsection{Preliminaries.}
First some general facts concerning $p(t)$-homogeneous operators are
given. Let an operator $A$ in $\sH$ be $p(t)$-homogeneous with
respect to a one-parameter family ${\mathfrak
U}=\{U_t\}_{t\in{\mathfrak{T}}}$ of unitary operators acting on
${\mathfrak H}$, cf. Definition \ref{dad1}. It follows from
(\ref{e1}) and (\ref{ee6}) that
\begin{equation}\label{tat1}
  p(t)p(g(t))=1, \qquad \forall{t}\in\mathfrak{T},
 \end{equation}
where the function of conjugation $g(t):
{\mathfrak{T}}\to{\mathfrak{T}}$ is determined by the formula
\begin{equation}\label{ee7}
  U_{g(t)}=U^*_t, \qquad \forall{t}\in\mathfrak{T}.
\end{equation}
\begin{lemma}\label{prop33}
Let $A$ be a $p(t)$-homogeneous operator with respect to a family
$\sU=\{U_t\}_{t\in{\mathfrak{T}}}$. Then for all
$t\in{\mathfrak{T}}$ and all $z\in\dC$,
\begin{equation}\label{eq0}
 U_t (\ker(A-zI))=\ker\left(p(t)A-zI\right).
\end{equation}
In particular, $\ker A$ is a reducing subspace for every $U_t$,
$t\in \sT$. Furthermore, $z\in\sigma_{a}(A) \iff \
z{p(t)}^n\in\sigma_{a}(A), \ n\in\mathbb{Z}, \ t\in\sT, \ a\in\{p,
r, c\}.$

If $p(t)\not=1$ at least for one point $t\in\sT$, then the essential
spectrum of $A$ contains the point $z=0$.
\end{lemma}
\begin{proof}
In view of (\ref{tat1}), $p(t)\not=0$ for all $t\in\sT$. Using
(\ref{ee6}) one gets
\begin{equation}\label{sas105}
 U_{t}({A}-z{I})=(p(t)A-zI)U_{t}
 =p(t)\left({A}-\frac{z}{p(t)}{I}\right)U_{t}
\end{equation}
that gives $U_t (\ker(A-zI))\subset \ker\left(p(t)A-zI\right)$. The
reverse inclusion is obtained by using \eqref{tat1}. The property of
$\ker A$ to be a reducing subspace for every $U_t$ follows from
(\ref{eq0}) with $z=0$ if one takes into account that $p(t)\not=0$.

The remaining assertions of the lemma immediately follow from
(\ref{sas105}).
\end{proof}
\begin{lemma}\label{x1}
Let $A$ be a closed densely defined $p(t)$-homogeneous operator with
respect to a family $\sU=\{U_t\}_{t\in{\mathfrak{T}}}$. Then also
its adjoint $A^*$ is
$p(t)$-homogeneous with respect to $\sU$. % and moreover for all
\end{lemma}
\begin{proof}
Since $A$ is $p(t)$-homogeneous one has $U_tA=p(t)AU_t$ for all
$t\in\sT$. As a unitary operator $U_t$ is bounded with bounded
inverse, and therefore, the previous equality is equivalent to
$
 A^*U_t^*=p(t)U_t^*A^* \quad \Longleftrightarrow \quad
 U_tA^*=p(t)A^*U_t$, \quad $\forall t\in\sT,$
which means that $A^*$ is $p(t)$-homogeneous with respect to $\sU$.
\end{proof}

In the case that $A$ is symmetric the formula (\ref{eq0}) in
Lemma~\ref{prop33} shows how the unitary operators $U_t$, $t\in\sT$,
transform the defect subspaces $\ker(A^*-zI)$ of $A$.

\begin{corollary}\label{x2}
Let $A$ in Lemma~\ref{x1} be nonnegative and $p(t)$-homogeneous with
respect to $\sU=\{U_t\}_{t\in{\mathfrak{T}}}$ and let $A_0$ be a
nonnegative selfadjoint extension of $A$. Then $
(p(t)A_0+I)(A_0+I)^{-1}U_t(\ker(A^*+I))=\ker(A^*+I).$
\end{corollary}
\begin{proof}
By Lemma~\ref{x1} the adjoint $A^*$ of $A$ is also
$p(t)$-homogeneous and \eqref{eq0} implies that $U_t(\ker(A^*+I))=
 \ker\left(A^*+1/{p(t)}\,I\right)$.
Moreover, the equality
$$
 (p(t)A_0+I)(A_0+I)^{-1}\ker\left(A^*+\frac{1}{p(t)}\,I\right)=\ker(A^*+I)
$$
is always satisfied for a nonnegative self-adjoint extension $A_0$
of $A$.
\end{proof}

For the next result recall that if $A$ is a nonnegative operator (or
in general a nonnegative relation) in a Hilbert space $\sH$, then
the Friedrichs extension $A_F$ and the Krein-von Neumann extension
$A_N$ of $A$ can be characterized as follows (see \cite{AN} for the
densely defined case and \cite{H, HMS, HSSW} for the general case):

If $\{ f , f '\} \in A^{*}$, then $\{ f, f'\} \in A_{F}$ if and only
if
\begin{equation}\label{FR}
\inf \left \{ \| f - h \|^{2} + (f'- h' , f - h) \, : \, \{h,h'\}
\in A \right\} =0.
\end{equation}

If $\{ f , f '\} \in A^{*}$, then $\{ f, f'\} \in A_{N}$ if and only
if
\begin{equation}\label{VN}
\inf \left \{ \| f' - h' \|^{2} + (f'- h' , f - h) \, : \, \{h,h'\}
\in A \right\} =0.
\end{equation}

\begin{lemma}\label{l34}
Let $A$ be a nonnegative densely defined $p(t)$-homogeneous operator
with respect to ${\mathfrak U}$. Then the Friedrichs extension $A_F$
and the Krein-von Neumann extension $A_N$ of $A$ are also
$p(t)$-homogeneous with respect to ${\mathfrak U}$. Moreover,
$U_t(\cD(A_F^{1/2}))\subset\cD(A_F^{1/2})$ and
$U_t(\cR(A_N^{1/2}))\subset\cR(A_N^{1/2})$ for all
${t}\in\mathfrak{T}$.
 \end{lemma}
\begin{proof} By Lemma~\ref{x1} $A^*$
is $p(t)$-homogeneous with respect to $\sU$. Hence, in view of
(\ref{e1}) and (\ref{ee6}), an intermediate extension $\wt A$ of $A$
is $p(t)$-homogeneous with respect to $\sU$ if and only if
\begin{equation}\label{les102}
 U_t : {\mathcal D}(\wt A) \to {\mathcal D}(\wt A), \quad
 \forall{t}\in\mathfrak{T}.
\end{equation}

To prove that $A_F$ is $p(t)$-homogeneous with respect to $\sU$,
assume that $f\in\cD(A_F)$. Then $g=U_tf\in \cD(A^*)$ and there is a
sequence $h_n\in\mathcal{D}(A)$ attaining the infimum in (\ref{FR}).
Then $U_th_n\in\mathcal{D}(A)$, $U_th_n\to U_tf=g$, and
\begin{equation}\label{eq3}
 \left( A^*U_tf-AU_th_n,U_tf-U_th_n\right)
 =(p(g(t)) \left( A^*f-Ah_n,f-h_n\right)\to 0,
\end{equation}
so that $g\in \cD(A_F)$ by (\ref{FR}). Therefore,
$U_t(\cD(A_F))\subset \cD(A_F)$ and $A_F$ is $p(t)$-homogeneous with
respect to $\sU$.

To prove the $p(t)$-homogeneity of $A_N$ assume that $f\in\cD(A_N)$.
Then again $g=U_tf\in\cD(A^*)$ and there is a sequence
$h_n\in\mathcal{D}(A)$ attaining the infimum in (\ref{VN}). In
particular, $Ah_n\to A^*f$, $U_th_n\in\mathcal{D}(A)$, and
$$
 AU_th_n=p(g(t))U_t Ah_n
 \to p(g(t))U_t A^*f=A*U_tf=A^*g.
$$
Moreover, (\ref{eq3}) is satisfied. Therefore, (\ref{VN}) shows that
$g\in \cD(A_N)$. This proves that $U_t(\cD(A_N))\subset \cD(A_N)$
and thus $A_N$ is $p(t)$-homogeneous with respect to $\sU$.

Finally, recall that the domain $\cD=\cD(A_F^{1/2})$, see
(\ref{kk2}), can be characterized as the set of vectors $f\in\sH$
satisfying
\[
 h_{n} \rightarrow f, \quad \left( A(h_{n}-h_{m}),h_{n}-h_{m}
\right) \rightarrow 0, \quad m, \, n \, \rightarrow \infty,
\]
and the range $\cR(A_N^{1/2})$ as the set of vectors $g\in\sH$
satisfying
\[
 Ah_{n} \rightarrow g,
 \quad \left(A(h_{n}-h_{m}),h_{n}-h_{m}
\right) \rightarrow 0, \quad m, \, n \, \rightarrow \infty,
\]
with $h_n\in \cD(A)$. The last statement is clear from these
characterizations using similar arguments as above with the sequence
$h_n$. This completes the proof.
\end{proof}

Let the operator $A_0$ in (\ref{ne3}) be $p(t)$-homogeneous with
respect to  ${\mathfrak U}=\{U_t\}_{t\in{\mathfrak{T}}}$. Define a
family of self-adjoint operators on $\sH$ by
\begin{equation}\label{tat11}
G_t=(p(t)A_0+I)(A_0+I)^{-1}, \quad {t}\in\mathfrak{T}.
\end{equation}
Clearly, $G_t$ is positive and bounded with bounded inverse for all
${t}\in\mathfrak{T}$. Moreover, it follows from (\ref{ee6}) and
(\ref{tat1}) that $(A_0+I)^{-1}U_t=U_t(p(g(t))A_0+I)^{-1}$  and
\begin{equation}\label{tat15}
 G_tU_t=U_tG^{-1}_{g(t)}=(G_{g(t)}U_{g(t)})^{-1}.
\end{equation}
Since $\|u\|_{-2}=\|(A_0+I)^{-1}u\|$, the identity
$(A_0+I)^{-1}U_t=G_tU_t(A_0+I)^{-1}$ implies that $
 \| U_t u \|_{-2} \le \|G_t\|\, \| u \|_{-2}$ for
all $u\in\mathfrak{H}$. Hence, the operators $U_t$ can be
continuously extended to bounded operators ${\mathbb{U}}_t$ in
${\mathfrak H}_{-2}(A_0)$ and, furthermore,
\begin{equation}\label{tat16}
 (\mathbb{A}_0+I)^{-1}{\mathbb{U}_t}\psi
 =G_tU_t(\mathbb{A}_0+I)^{-1}\psi
\end{equation}
for all $\psi\in{\mathfrak H}_{-2}(A_0)$ and $t\in\sT$. The equality
(\ref{ee7}) shows that $\dU_t$ has a bounded inverse which satisfies
$\dU_t^{-1}=\dU_{g(t)}$. The operator $\dU_{t}$ can be characterized
also as the dual mapping (adjoint) of $U_{g(t)}$ with respect to the
form defined in (\ref{ada5}). In fact, using (\ref{ada5}),
(\ref{ee6}), (\ref{ee7}), and (\ref{tat16}), it is seen that the
action of the functional $<{\mathbb{U}}_t\psi, \cdot>$ on the
elements $u\in{\mathfrak H}_{2}(A_0)$ is determined by the formula
\begin{eqnarray}\label{tat2}
 <\mathbb{U}_t\psi, u>=((A_0+I)u, G_tU_th)=(U_{g(t)}(p(t)A_0+I)u, h) & & \nonumber \\
 =((A_0+I)U_{g(t)}u, h)=<\psi, U_{g(t)}u>,
\end{eqnarray}
where $h=(\mathbb{A}_0+I)^{-1}\psi$.

Now consider a singular element $\psi\in \sH_{-2}(A_0)$, cf.
(\ref{ne3}). The assumption that $\psi$ is $\xi(t)$-invariant with
respect to ${\mathfrak U}$, i.e. $\dU_t\psi=\xi(t)\psi$ for all
$t\in\sT$ (see Definition \ref{dad2}), implies some relations
between $\xi(t)$, $p(t)$, and $g(t)$.

\begin{proposition}\label{p1}
Let the operator $A_0$ in (\ref{ne3}) be $p(t)$-homogeneous with
respect to the family $\sU$ and let $\psi\in{\mathfrak
H}_{-2}(A_0)\setminus{\mathfrak  H}$ be $\xi(t)$-invariant with
respect to $\sU$. Then for all ${t}\in\mathfrak{T}$ one has
\begin{equation}\label{tat1b}
 \xi(t)\xi(g(t))=1
\end{equation}
and, moreover, $|\xi(t)|=1$ if $p(t)=1$ and $\min \{1,
p(t)\}<|\xi(t)|<\max \{1,p(t)\}$ if $p(t)\not=1$.
\end{proposition}
\begin{proof}
It follows from (\ref{tat10}) and (\ref{tat16}) that
$\psi\in{\mathfrak  H}_{-2}(A_0)\setminus{\mathfrak H}$ is
$\xi(t)$-invariant with respect to ${\mathfrak U}$ if and only if
\begin{equation}\label{tat12}
 G_tU_th=\xi(t)h,\quad \forall{t}\in\mathfrak{T},
\end{equation}
where $h=(\mathbb{A}_0+I)^{-1}\psi$. This together with
(\ref{tat15}) implies that
$$
 h=(G_{g(t)}U_{g(t)})(G_tU_t) h
  =\xi(t)G_{g(t)}U_{g(t)}h=\xi(t)\xi(g(t))h,
$$
which proves (\ref{tat1b}). Moreover, (\ref{tat12}) shows that
$|\xi(t)|\|h\|=\|G_tU_th\|.$ In particular, if $p(t)=1$, then
$G_t=I$ and $|\xi(t)|\|h\|=\|U_th\|=\|h\|$ that gives $|\xi(t)|=1$.

In the case where $p(t)\not=1$ the formula for $G_t$ in
(\ref{tat11}) with an evident reasoning leads to the estimates
$$
 \alpha(t)\|h\|=\alpha(t)\|U_th\|<\|G_tU_th\|<\beta(t)\|U_th\|=\beta(t)\|h\|,
$$
where $\alpha(t)=\min\{1, p(t)\}$ and $\beta(t)=\max\{1,p(t)\}$.
This completes the proof.
\end{proof}

\subsection{$p(t)$-homogeneous self-adjoint extensions of $A_{\mathrm{sym}}$.}
Let $A_{\mathrm{sym}}$ be defined by (\ref{e7}). This means that
$A_{\mathrm{sym}}$ is a nonnegative symmetric operator with finite
defect numbers.
\begin{lemma}\label{t33}
If $p(t)\not=1$ at least for one point $t\in\sT$, then an arbitrary
$p(t)$-homogeneous self-adjoint extension of the symmetric operator
$A_{\mathrm{sym}}$ is nonnegative.
\end{lemma}
\begin{proof}
Assume that $z$ is a negative eigenvalue of a $p(t)$-homogeneous
self-adjoint extension $A$ of $A_{\mathrm{sym}}$ and that
$p(t)\not=1$ for $t\in{\mathfrak{T}}$.  Then, according to Lemma
\ref{prop33}, there exists infinite series of negative eigenvalues
$z{p(t)^n}$ \ $(n\in{\mathbb Z})$ of $A$ that contradicts to the
assumption of finite defect numbers of $A_{\mathrm{sym}}$. Hence,
$A$ is a nonnegative extension of $A_{\mathrm{sym}}$.
\end{proof}
\begin{lemma}\label{l56}
Let $A_0$ be $p(t)$-homogeneous and let $\psi_j$ be
$\xi_j(t)$-invariant with respect to $\sU$, $j=1,\dots,n$. Then the
symmetric operator $A_{\mathrm{sym}}$ defined by (\ref{e7}) and its
adjoint $A_{\mathrm{sym}}^*$ are also $p(t)$-homogeneous with
respect to ${\mathfrak U}$.
\end{lemma}
\begin{proof}
It follows from (\ref{e7}) and (\ref{tat2}) that
  $$
 <\psi_j, U_tu>=<\mathbb{U}_{g(t)}\psi_j, u>=\xi_j(g(t))<\psi_j, u>=0
  $$
for every $u\in{\mathcal D}(A_{\mathrm{sym}})$. Thus
 $U_t:{\mathcal D}(A_{\mathrm{sym}})\to{\mathcal
 D}(A_{\mathrm{sym}})$ and hence by (\ref{ee6})
$A_{\mathrm{sym}}$ is $p(t)$-homogeneous:
$U_tA_{\mathrm{sym}}=p(t)A_{\mathrm{sym}}U_t$. By Lemma~\ref{x1}
also the adjoint $A_{\mathrm{sym}}^*$ is $p(t)$-homogeneous with
respect to $\sU$.
\end{proof}
%If the assumptions in Lemma~\ref{l56} are satisfied, the defect
%subspace $\ker(A_{\mathrm{sym}}^*+I)$ of $A_{\mathrm{sym}}$ is
%invariant under $G_tU_t$, see Corollary~\ref{x2}.

In view of (\ref{tat10}) and (\ref{tat2}) the $\xi_j(t)$-invariance
of $\psi_j$ is equivalent to the relation
\begin{equation}\label{tat3}
 \xi_j(t)<\psi_j, u>=<\psi_j, U_{g(t)}u>, \quad
 \forall{u}\in{\mathfrak H}_{2}(A_0), \quad
 \forall{t}\in\mathfrak{T},
\end{equation}
where the linear functionals $<\psi_j,\cdot>$ are defined by
(\ref{ada5}). The next theorem shows that the preservation of
(\ref{tat3}) for the extended functionals $<\psi_j^{\mathrm{ex}},
\cdot>$  is closely related to the existence of $p(t)$-homogeneous
self-adjoint extensions of $A_{\mathrm{sym}}$ transversal to $A_0$.

\begin{theorem}\label{t5}
Let $A_0$ be $p(t)$-homogeneous, let $\psi_1,\dots,\psi_n$ be
$\xi_j(t)$-invariant with respect to $\sU$, and let
$<\psi_j^{\mathrm{ex}}, f>$  be defined by (\ref{k23}). Then the
relations
\begin{equation}\label{les100}
 \xi_j(t)<\psi_j^{\mathrm{ex}}, f>=<\psi_j^{\mathrm{ex}}, U_{g(t)}f>,
 \quad 1\leq{j}\leq{n},  \quad  \forall{t}\in\mathfrak{T},
\end{equation}
are satisfied for all $f\in{\mathcal D}(A_{\mathrm{sym}}^*)$ if and
only if the corresponding self-adjoint operator $\widetilde{A}$
defined by \eqref{lesia200} is $p(t)$-homogeneous with respect to
${\mathfrak U}$.
\end{theorem}
\begin{proof}
Denote
\begin{equation}\label{did4}
{\mathbf{\Xi}}(t)=\left(\begin{array}{cccc}
\xi_1(t) & 0 & \ldots & 0 \\
0 & \xi_2(t) & \ldots & 0 \\
\vdots & \vdots & \ddots & \vdots \\
0 & 0 & \ldots & \xi_n(t)
\end{array}\right).
\end{equation}
Then $\det\mathbf{\Xi}(t)\neq 0$, $t\in\sT$, by
Proposition~\ref{p1}, since $\psi_i$ is $\xi_j(t)$-invariant with
respect to $\sU$. By using (\ref{k9}) in Lemma~\ref{l23} the
relations (\ref{les100}) can be rewritten as follows:
\begin{equation}\label{les101}
 {\mathbf{\Xi}}(t)\Gamma_0f=\Gamma_0U_{g(t)}f,
 \quad \forall{f}\in{\mathcal D}(A_{\mathrm{sym}}^*),
 \quad \forall{t}\in\mathfrak{T}.
\end{equation}
Since ${\mathcal D}(\widetilde{A})=\ker\Gamma_0$, (\ref{les101})
immediately implies that $U_t(\cD(\wt A))\subset\cD(\wt A)$, cf.
(\ref{ee7}). Thus the equalities (\ref{les100}) ensure the
$p(t)$-homogeneity of $\widetilde{A}$ with respect to ${\mathfrak
U}$.

Conversely, assume that $\wt A$ is $p(t)$-homogeneous with respect
to ${\mathfrak U}$. According to (\ref{lesia200}), (\ref{ee7}), and
(\ref{les102}) this is equivalent to
\begin{equation}\label{les103}
 -{\mathbf R}\widehat{\Gamma}_0U_{g(t)}f=\widehat{\Gamma}_1U_{g(t)}f,
 \quad \forall{f}\in{{\mathcal D}(\widetilde{A})}, \quad
 \forall{t}\in\mathfrak{T}.
\end{equation}
Using (\ref{tat11}), (\ref{tat1b}), and (\ref{tat12}) it is seen
that
\begin{equation}\label{tat55}
\begin{split}
 U_{g(t)}h_j
 &=p(t)G_{g(t)}U_{g(t)}h_j+(I-p(t)G_{g(t)})U_{g(t)}h_j \\
 &=\frac{p(t)}{\xi_j(t)}h_j+(1-p(t))(A_0+I)^{-1}U_{g(t)}h_j,
 \end{split}
\end{equation}
where $h_j=(\mathbb{A}_0+I)^{-1}\psi_j$, $j=1,\ldots,n$. This
expression and relations (\ref{lesia99}), (\ref{tat2}) yield the
following equalities for all $f=u+\sum_{j=1}\alpha_jh_j\in{\mathcal
D}(A_{\mathrm{sym}}^*)$ and ${t}\in\mathfrak{T}$:
\begin{equation}\label{les10}
 \widehat{\Gamma}_0U_{g(t)}f=p(t){\mathbf{\Xi}}^{-1}(t)\widehat{\Gamma}_0f,
 \quad \widehat{\Gamma}_1U_{g(t)}f
  =\mathbf{\Xi}(t)\widehat{\Gamma}_1f+(1-p(t)){\mathbf{G}}^{\top}(t)
   \widehat{\Gamma}_0f,
\end{equation}
where ${\mathbf{G}}^{\top}(t)$ is the transpose of the matrix
${\mathbf{G}}(t)=((h_i, U_th_j))_{i,j=1}^n$. Now with $f\in{\mathcal
D}(\widetilde{A})$ substituting these expressions into
(\ref{les103}), using (\ref{lesia200}), and taking into account that
$\widehat{\Gamma}_0(\cD(\widetilde{A}))=\dC^n$, one concludes that
the $p(t)$-homogeneity of $\widetilde{A}$ is equivalent to the
matrix equality
\begin{equation}\label{les104}
 {\mathbf{\Xi}}(t){\mathbf{R}}-p(t){\mathbf{R}}{\mathbf{\Xi}}^{-1}(t)
 =(1-p(t)){\mathbf{G}}^{\top}(t),
 \quad \forall{t}\in{\mathfrak{T}}.
\end{equation}
Finally, employing (\ref{lesia101}) and (\ref{les10}) it is easy to
see that equality (\ref{les104}) is equivalent to (\ref{les101}).
Therefore, the extended functionals $<\psi_j^{\mathrm{ex}}, \cdot>$
satisfy the relations (\ref{les100}). Theorem \ref{t5} is proved.
\end{proof}
\begin{remark}
In the particular case where $p(t)=t^\beta$ and $\xi(t)=t^{\theta}$
with $\beta, \theta\in{\mathbb{R}}$, another condition for the
preservation of $\xi(t)$-invariance for $<\psi_j^{\mathrm{ex}},
\cdot>$ has been obtained in \cite[Lemma 1.3.2]{AL1}.
\end{remark}
\begin{corollary}\label{new28}
Let $\widetilde{A}$ be a self-adjoint extension of
$A_{\mathrm{sym}}$ transversal to $A_0$. Then $\widetilde{A}$ is
$p(t)$-homogeneous if and only if $\widetilde{A}$ is defined by
(\ref{lesia200}) and the entries $r_{ij}$ of $\mathbf{R}$ in
(\ref{lesia200}) satisfy the following system of equations for all
${t}\in{\mathfrak{T}}$:
\begin{equation}\label{kak9}
 \beta_{ij}(t)r_{ij}=(1-p(t))(h_j, U_th_i),
 \quad
 \beta_{ij}(t)=\left(\xi_i(t)-\frac{p(t)}{\xi_j(t)}\right),
 \quad
 1\leq{i,j}\leq{n}.
\end{equation}
\end{corollary}
\begin{proof} Since $\ker\widehat{\Gamma}_0=\cD(A_0)$, formula (\ref{lesia200})
describe all self-adjoint extensions of $A_{\mathrm{sym}}$
transversal to $A_0$ when the parameter
$\mathbf{R}=(r_{ij})_{i,j=1}^n$ runs the set of all Hermitian
matrices. Hence, $\widetilde{A}=A_\mathbf{R}$ for some choice of
$\mathbf{R}$ in (\ref{lesia200}). The proof of Theorem~\ref{t5}
shows that $A_\mathbf{R}$ is $p(t)$-homogeneous if and only if
$\mathbf{R}$ is a solution of (\ref{les104}) that does not depend on
${t}\in{\mathfrak{T}}$. Rewriting (\ref{les104}) componentwise one
gets (\ref{kak9}).
\end{proof}
\begin{remark}\label{new44}
In the case that $p(x)\equiv{1}$, the right-hand side of
(\ref{kak9}) vanishes and (\ref{kak9}) reduces to
$\beta_{ij}(t)r_{ij}=0$, $1\leq{i,j}\leq{n}$. Moreover, by
Proposition~\ref{p1} $\beta_{ii}(t)\equiv{0}$ and, therefore, the
entries $r_{ii}$ cannot be uniquely determined from (\ref{kak9}).
This implies the existence of infinitely many $1$-homogeneous
self-adjoint extensions of $A_{\mathrm{sym}}$ transversal to $A_0$.
\end{remark}

\begin{example}\label{examp1}
Let $\alpha>0$ and let $\wt A$ be defined by
$$
 \widetilde{A}_\alpha=A_{\mathrm{sym}}^*\upharpoonright
  {\mathcal{D}(\widetilde{A}_\alpha)},
 \quad
 \mathcal{D}(\widetilde{A}_\alpha)
  =\mathcal{D}(A_{\mathrm{sym}})\dot{+}\ker(A_{\mathrm{sym}}^*+\alpha{I}).
$$
Then for all $\alpha>0$, $\widetilde{A}_\alpha$ is a $1$-homogeneous
self-adjoint extensions of $A_{\mathrm{sym}}$ transversal to $A_0$.
\end{example}

\subsection{Uniqueness of $p(t)$-homogeneous admissible operators.}
Let the operator $A_0$ be $p(t)$-homogeneous and let the singular
elements $\psi_j$ appearing in (\ref{ne3}) be $\xi_j(t)$-invariant
with respect to ${\mathfrak U}$.

If all $\psi_j$ belong to ${\mathfrak H}_{-1}(A_0)$, then the
extended functionals $<\psi_j^{\mathrm{ex}}, \cdot>$ are determined
by continuity onto $\mathcal{D}(A_{\mathrm{sym}}^*)$ and they
automatically possess the property of $\xi_j(t)$-invariance
(\ref{les100}), since $U_t\upharpoonright_{\mathcal{D}(A_0)}$ can be
extended by continuity onto  ${\mathfrak H}_{1}(A_0)$. In this case,
the set of admissible operators consists of a unique element (the
Friedrichs extension $A_F$, see Corollary \ref{p121}) and this
admissible operator is $p(t)$-homogeneous.

If  ${\mathfrak H}_{-1}(A_0)$ does not contain all $\psi_j$, then
admissible operators for the regularization of (\ref{ne3}) are not
determined uniquely. In this case, the natural assumption of
$\xi_j(t)$-invariance for the extended functionals
$<\psi_j^{\mathrm{ex}}, \cdot>$ can be used to select a unique
admissible operator $\widetilde{A}$. By Theorem~\ref{t5} the
$\xi_j(t)$-invariance of $<\psi_j^{\mathrm{ex}}, \cdot>$ is
equivalent to the $p(t)$-homogeneity of the corresponding operator
$\widetilde{A}$ defined by (\ref{lesia200}). Therefore, instead of
assumption of $\xi_j(t)$-invariance one can use the requirement of
$p(t)$-homogeneity imposed on the set of admissible operators to
achieve their uniqueness.

\begin{theorem}\label{new2007}
Assume that the singular elements $\psi_j$ in (\ref{ne3}) are
${\mathfrak H}_{-1}(A_0)$-independent and the system of equations
(\ref{kak9}) has a unique solution $\mathbf{R}=(r_{ij})_{i,j=1}^n$
that does not depend on $t\in\mathfrak{T}$. Then there exists a
unique $p(t)$-homogeneous admissible operator $\widetilde{A}$ for
the regularization of (\ref{ne3}) and it coincides with the
Krein-von Neumann extension $A_N$ of $A_{\mathrm{sym}}$.
 \end{theorem}
\begin{proof}
Let $\mathbf{R}=(r_{ij})_{i,j=1}^n$ be a unique solution of
(\ref{kak9}) and let $\widetilde{A}$ be the corresponding
self-adjoint extension of $A_{\mathrm{sym}}$ determined by
(\ref{lesia200}).

Since (\ref{kak9}) has a unique solution, $p(t)\not=1$ for at least
one point $t\in\mathfrak{T}$ (see Remark \ref{new44}). In this case,
Lemma \ref{t33} and relation (\ref{lesia200}) imply that
$\widetilde{A}$ is a nonnegative extension of $A_{\mathrm{sym}}$
transversal to $A_0$. Then also $A_F$ and $A_N$ are transversal
extensions of $A_{\mathrm{sym}}$; cf. the proof of Theorem
\ref{t303}. These extensions are also $p(t)$-homogeneous (see Lemmas
\ref{l56}, \ref{l34}).

Since elements $\psi_j$ in \eqref{ne3} form an ${\mathfrak
H}_{-1}(A_0)$-independent system, Corollary \ref{p1211} gives that
any self-adjoint extension of $A_{\mathrm{sym}}$ transversal to
$A_0$ is admissible for the regularization of (\ref{ne3}) and
$A_0=A_F$. The unique solution of (\ref{kak9}) allows one to select
a unique $p(t)$-homogeneous self-adjoint extension $\widetilde{A}$
of $A_{\mathrm{sym}}$ transversal to $A_0=A_F$. Obviously, it
coincides with the Krein-von Neumann extension $A_N$.
\end{proof}

The next statement concerns to the general case.
\begin{theorem}\label{did11}
Let $A_F$ and $A_N$ be transversal, let the operator $S$ defined in
(\ref{les36}) be  $p(t)$-homogeneous for some choice of
$\widetilde{\mathcal{H}}$ satisfying conditions (\ref{dur57}), and
assume that for every $\beta_{ij}(t)$ in (\ref{kak9}) there exists
at least one point $t_{ij}\in{\mathfrak{T}}$ such that
$\beta_{ij}(t_{ij})\not=0$. Then there exists a unique
$p(t)$-homogeneous admissible operator for the regularization of
(\ref{ne3}).
 \end{theorem}
\begin{proof}
Let $\widetilde{A}$ be the Krein-von Neumann extension of $S$. The
second part of the proof of Theorem \ref{t303} shows that
$\widetilde{A}$ is an admissible operator. By Lemma \ref{l34},
$\widetilde{A}$ is $p(t)$-homogeneous. Its uniqueness follows from
the fact that condition $\beta_{ij}(t_{ij})\not=0$ ensures in view
of (\ref{kak9}) the uniqueness of $p(t)$-homogeneous self-adjoint
extensions of $A_{\mathrm{sym}}$ transversal to $A_0$.
\end{proof}

The next statement contains conditions for the $p(t)$-homogeneity of
the symmetric operator $S$ defined by (\ref{les36}) in Lemma
\ref{lil1} which appear to be useful in applications.

\begin{proposition}\label{pip2}
Let $A_0$ be $p(t)$-homogeneous, let the singular elements $\psi_j$
in (\ref{ne3}) be $\xi_j(t)$-invariant with respect to $\sU$, and
let $\cY=(\dA_0+I)(\mathcal{H}\ominus\widetilde{\mathcal H})$. Then:
\begin{enumerate}
\def\labelenumi{\rm (\roman{enumi})}

\item $S$ is $p(t)$-homogeneous if and only if
$\cY$ is invariant under $\dU_t$, $t\in \sT$, and
\begin{equation}\label{sas60}
(h', U_t\wt{h}^{\perp})=0, \ \  \forall{h'}\in{\mathcal{H}}',  \ \
 \forall{\wt{h}^{\perp}}\in\mathcal{H}\ominus\widetilde{\mathcal H}, \ \
 \forall{t}\in\mathfrak{T}_0=\{\, t\in\mathfrak{T} : \ p(t)\not=1 \}.
\end{equation}

\item If $G_tU_t$, $t\in \sT$, is self-adjoint, then
$S$ with $\widetilde{\mathcal H}={\mathcal H}'$ is
$p(t)$-homogeneous if and only if (\ref{sas60}) holds.

\item If $\cY$ is a linear span of some singular
elements $\psi_j$ in (\ref{ne3}), then $S$ is $p(t)$-homogeneous if
and only if (\ref{sas60}) holds.
\end{enumerate}
\end{proposition}

\begin{proof}
(i) The definition (\ref{les36}) shows that
$\ker(S^*+I)=\cH\ominus\wt\cH$. Hence, if $S$ is $p(t)$-homogeneous
with respect to ${\mathfrak U}$ then
$G_tU_t(\cH\ominus\wt\cH)=\cH\ominus\wt\cH$ by Corollary~\ref{x2}.
According to (\ref{tat16}) the subspace $\cH\ominus\wt\cH$ is
invariant under $G_tU_t$ if and only if
$\cY=(\dA_0+I)(\cH\ominus\wt\cH)$ is invariant under the operator
$\dU_t$, $t\in\sT$. Thus, if $S$ is $p(t)$-homogeneous with respect
to ${\mathfrak U}$ then $\cY$ is invariant under $\dU_t$, $t\in\sT$.

By Lemma \ref{l56}, $A_{\mathrm{sym}}^*$ is $p(t)$-homogeneous with
respect to ${\mathfrak U}$. Since $S$ is an intermediate extension
of $A_{\mathrm{sym}}$ its $p(t)$-homogeneity is equivalent to the
relation $U_{g(t)}(\mathcal{D}(S))\subset\mathcal{D}(S),$
${t}\in\mathfrak{T}$, see (\ref{les102}).

The definition of $S$ in (\ref{les36}) implies that
\begin{equation}\label{sas71}
 U_{g(t)}f\in{\mathcal D}(S)\iff((A_F+I)U_{g(t)}f, \wt{h}^{\perp})=0,
  \quad \forall\wt{h}^{\perp}\in\mathcal{H}\ominus\widetilde{\mathcal{H}}.
\end{equation}

Now let $f=h'+u\in{\mathcal D}(S)$ be decomposed as in
Lemma~\ref{lil1}, see (\ref{sas90}), (\ref{sas70}). It follows from
(\ref{tat55}) that
$$
(A_F+I)U_{g(t)}f=(A_{\mathrm{sym}}^*+I)U_{g(t)}f=(1-p(t))U_{g(t)}h'+(A_0+I)U_{g(t)}u.
$$
By taking (\ref{tat2}) into account one obtains
\begin{equation}\label{eq1}
\begin{split}
 ((A_F+I)U_{g(t)}f, \wt{h}^{\perp})
  &=(1-p(t))(U_{g(t)}h', \wt{h}^{\perp})+((A_0+I)U_{g(t)}u, \wt{h}^{\perp}) \\
  &=(1-p(t))(h', U_{t}\wt{h}^{\perp})+ < \dU_t \psi, u>.
\end{split}
\end{equation}
If $\cY$ is invariant under $\dU_t$, $t\in\sT$, then $<\dU_t \psi,
u>=0$ for all $f=h'+u\in{\cD}(S)$. Now (\ref{sas71}) and (\ref{eq1})
show that $S$ is $p(t)$-homogeneous if and only if $\cY$ is
invariant under $\dU_t$ and (\ref{sas60}) holds.

(ii) Since $A_0$ and $A_F$ are $p(t)$-homogeneous, the symmetric
restriction $S_0:=A_F\cap A_0$ and its adjoint $S_0^*$ are also
$p(t)$-homogeneous, see Lemma~\ref{x1}. It follows from (\ref{tut})
that $f\in\cD(S_0)$ if and only if $f\in\cD(A_0)$ and
\[
  ((A_0+I)f,h')=0,\quad \forall h'\in\cH'=\cH\cap \sH_1(A_0).
\]
Hence, $\ker(S_0^*+I)=\cH'$ and $G_tU_t\cH'=\cH'$ for all $t\in\sT$
by Corollary~\ref{x2}. Similarly $G_tU_t\cH=\cH$ for all $t\in\sT$,
since $A_{\mathrm{sym}}$ is $p(t)$-homogeneous. Therefore, if
$G_tU_t$ is self-adjoint, then $\cH$ and $\cH'$ are reducing
subspaces for the operators $G_tU_t$ and consequently
$G_tU_t\cH''\subset\cH''$ is satisfied for all $t\in\sT$. Then,
according to (\ref{tat16}), $\cY=(\dA_0+I)\cH''$ is invariant under
$\dU_t$. Now the claim follows from part (i) with
$\widetilde{\mathcal H}={\mathcal H}'$ and
$\mathcal{H}\ominus\widetilde{\mathcal H}={\mathcal H}''$.

(iii) If $\cY$ has a basis formed by some $\xi_j(t)$-invariant
singular elements $\psi_j$, then $\cY$ is invariant under $\dU_t$,
see (\ref{tat10}). So, the statement is reduced to (i).
\end{proof}

\begin{example}\label{Example2}
{\it A general zero-range potential}. A one-dimensional
Schr\"{o}dinger operator corresponding to a general zero-range
potential at the point $x=0$ can be given by the expression
$$
  A_0+b_{11}<\delta,\cdot>\delta(x)+b_{12}<\delta',\cdot>\delta(x)+
  b_{21}<\delta,\cdot>\delta'(x)+b_{22}<\delta',\cdot>\delta'(x),
$$
where $A_0=-d^2/dx^2$ $(\mathcal{D}(A_0)=W_2^2(\mathbb{R}))$ acts in
${\mathfrak H}=L_2(\mathbb{R})$, $\delta'(x)$ is the derivative of the Dirac $\delta$-function
(with support at $0$).

In this case,
$A_{\mathrm{sym}}=-d^2/dx^2\upharpoonright{\{u(x)\in{W}_2^2(\mathbb{R})
:\, u(0)=u'(0)=0\}}$ and the corresponding Friedrichs and Krein-von
Neumann extensions are transversal (see, e.g., \cite{AHSS}). The
functions
$$
(\mathbb{A}_0+I)^{-1}\psi_1=h'(x)=\frac{1}{2}\left\{\begin{array}{cc}
 e^{-x}, & x>0  \\
 e^{x}, & x<0
 \end{array}\right. ,
$$
$$
  \ \ \
 (\mathbb{A}_0+I)^{-1}\psi_2=h''(x)=\frac{1}{2}\left\{\begin{array}{cc}
 -e^{-x}, & x>0  \\
 e^{x}, & x<0
 \end{array}\right. ,
$$
where $\psi_1=\delta(x)$ and $\psi_2=\delta'(x)$, form an orthogonal
basis of $\mathcal{H}=\ker(A_{\mathrm{sym}}^*+I)$ such that
${\mathcal{H}}^{\prime}=<h'(x)>$ and ${\mathcal{H}}^{''}=<h''(x)>$.

Define ${\mathfrak U}=\{U_t\}_{t\in{[0, \infty)}}$ as a collection
of \textit{the space parity operator} $U_0f(x)=f(-x)$
$(f(x)\in{L_2(\mathbb{R})})$ and the set of \textit{scaling
transformations} $U_tf(x)=\sqrt{t}f(tx)$,  $t>0$. In this case,
$A_0$ is $p(t)$-homogeneous with respect to ${\mathfrak U}$, where
$p(0)=1$ and $p(t)=t^{-2}$ if $t>0$.  The elements $\psi_j$ $(j=1,
2)$ are $\xi_j(t)$-invariant, where $\xi_1(0)=1$,
$\xi_1(t)=t^{-1/2}$ \ $(t>0)$
 and $\xi_2(0)=-1$, $\xi_2(t)=t^{-3/2}$ \ $(t>0)$.
Furthermore, for such a choice of ${\mathfrak U}$,
$\mathfrak{T}_0=\{\, t\in[0, \infty) : \ p(t)\not=1 \}=(0,\infty)$
and
$$
(h',
U_th'')=t^{1/2}\int_{-\infty}^{\infty}h'(x)\overline{h''(tx)}dx=0, \
\quad  \forall{t}\in\mathfrak{T}_0.
$$

Let us put $\widetilde{\mathcal{H}}={\mathcal H}'$. Then
$\cY=(\mathbb{A}_0+I){\mathcal{H}}^{''}=<\psi_2>$ and part (iii) of
Proposition~\ref{pip2} implies that the corresponding operator $S$
defined by (\ref{les36}) is $p(t)$-homogeneous. Calculating
$\beta_{ij}(t)$ in (\ref{kak9}) for $\xi_1(t)$, $\xi_2(t)$, and
$p(t)$ as given above, it is easy to see that $\beta_{ij}(0)\not=0$
if $i\not={j}$ and $\beta_{ii}(t)\not=0$ for all $t>0$. In this
case, by Theorem~\ref{did11} there exists a unique
$p(t)$-homogeneous admissible operator $\widetilde{A}$.

To identify $\widetilde{A}$ it suffices to determine the entries
$r_{ij}$ of the corresponding matrix $\mathbf{R}$ in
(\ref{lesia200}) with the aid of (\ref{kak9}):

For $t=0$, \eqref{kak9} takes the form $\left(\begin{array}{cc}
0 &  2r_{12} \\
-2r_{21} & 0
\end{array}\right)=0$ and, hence, $r_{12}=r_{21}=0$.
On the other hand, for $t>0$ calculating both sides of (\ref{kak9})
leads to
$$
t^{-3/2}(t-1)\left(\begin{array}{cc}
r_{11} &  0 \\
0 & -r_{22}
\end{array}\right)=(1-t^{-2})\left(\begin{array}{cc}
\frac{\sqrt{t}}{2(1+t)} &  0 \\
0 & \frac{\sqrt{t}}{2(1+t)}
\end{array}\right)$$
and thus $r_{11}=1/2, r_{22}=-1/2$. Substituting the coefficients
$r_{ij}$ in (\ref{k23}) results in the well-known extensions of
$\delta(x)$ and $\delta'(x)$ onto
$\mathcal{D}(A_{\mathrm{sym}}^*)=W_2^2(\mathbb{R}\backslash\{0\})$
(see \cite{AL1}):
$$
<\delta_{\mathrm{ex}}, f>=\frac{f(+0)+f(-0)}{2}, \hspace{5mm} <\delta_{\mathrm{ex}}', f>=-\frac{f'(+0)+f'(-0)}{2}.
$$
The corresponding $p(t)$-homogeneous admissible operator
$\widetilde{A}$ is the restriction of $-d^2/dx^2$ to $
 \mathcal{D}(\widetilde{A})=\left\{\,f(x)\in{W_2^2(\mathbb{R}\backslash\{0\})}
 :\,  -f(-0)=f(+0), \ -f'(-0)=f'(+0) \,\right\}.$
\end{example}

\subsection{The case of rank one singular perturbations.}
In the case of rank one singular perturbations
$A_0+b<\psi,\cdot>\psi$, where $A_0$ is $p(t)$-homogeneous and
$\psi$ is $\xi(t)$-invariant, the system (\ref{kak9}) takes the form
\begin{equation}\label{ada7}
(\xi^2(t)-{p(t)})r=\xi(t)(1-p(t))(h, U_th)
 \quad (h=(\mathbb{A}_0+I)^{-1}\psi),
\qquad \forall{t}\in\mathfrak{T}.
\end{equation}
\begin{proposition}\label{pip1}
1. If (\ref{ada7}) has no solutions, then there is only one
$p(t)$-homogeneous extension $A_0=A_F=A_N$ and any self-adjoint
extension of ${A_{\mathrm{sym}}}$ different from $A_0$ has a
negative eigenvalue.

2. If (\ref{ada7}) has at least two solutions, then all self-adjoint
extensions of $A_{\mathrm{sym}}$ are $p(t)$-homogeneous.

3. If (\ref{ada7}) has a unique solution $r\in\mathbb{R}$ that does
not depend on ${t}\in\mathfrak{T}$, then the symmetric operator
$A_{\mathrm{sym}}$ associated with $A_0+b<\psi,\cdot>\psi$ possesses
exactly two $p(t)$-homogeneous extensions: the Friedrichs $A_F$ and
the Krein-von Neumann $A_N$ extensions. One of them coincides with
$A_0$, another one is the unique $p(t)$-homogeneous admissible
operator $\widetilde{A}$ for the regularization of
$A_0+b<\psi,\cdot>\psi$. More precisely, $A_0=A_F$ and
$\widetilde{A}=A_N$ if $\psi\in{\mathfrak
H}_{-2}(A_0)\setminus{\mathfrak H}_{-1}(A_0)$; \ $A_0=A_N$ and
$\widetilde{A}=A_F$ if $\psi\in{\mathfrak H}_{-1}(A_0)$.
\end{proposition}
\begin{proof}
In the case of rank one perturbations, an arbitrary self-adjoint
extension ${A}(\neq{A_0})$ of the symmetric operator
$A_{\mathrm{sym}}=A_0\upharpoonright{\{\,u\in{\mathcal{D}}(A_0)
  :\, <\psi,u>=0 \,\}}$ is transversal to $A_0$. This means
that there is a one-to-one correspondence between the set of
solutions $r\in\mathbb{R}$ of (\ref{ada7}) and the set of
$p(t)$-homogeneous self-adjoint extensions $A(\neq{A_0})$ of
${A_{\mathrm{sym}}}$.

By Lemmas \ref{l34}, \ref{l56} the symmetric operator
$A_{\mathrm{sym}}$ and its Friedrichs $A_F$ and Krein-von Neumann
$A_N$ extensions are $p(t)$-homogeneous. Therefore, if (\ref{ada7})
has no solutions, then $A_N=A_F=A_0$ that justifies assertion 1.

Two different solutions of (\ref{ada7}) may appear only in the case
where $\xi^2(t)={p(t)}$ and $(1-p(t))(h, U_th)=0$ for all
$t\in\mathfrak{T}$. But these equalities are equivalent to the fact
that any $r\in\mathbb{R}$ is a solution of (\ref{ada7}). Therefore,
an arbitrary self-adjoint extension of $A_{\mathrm{sym}}$ is
$p(t)$-homogeneous. Assertion 2 is proved.

Finally, assume that (\ref{ada7}) has a unique solution. It follows
from Corollary \ref{new28} that the set of all
 $p(t)$-homogeneous extensions of ${A_{\mathrm{sym}}}$ is exhausted
 by the Friedrichs $A_F$ and the Krein-von Neumann $A_N$ extensions.
One of them coincides with $A_0$, another one is the unique
$p(t)$-homogeneous admissible operator $\widetilde{A}$. To complete
the proof it suffices to use Theorem \ref{new2007} for
$\psi\in{\mathfrak H}_{-2}(A_0)\setminus{\mathfrak H}_{-1}(A_0)$ and
Corollary \ref{p121} for $\psi\in{\mathfrak H}_{-1}(A_0)$.
\end{proof}

\begin{example}\label{Example1}
{\it One point interaction in} ${\mathbb R}^n$ $(n=1,2,3)$.

Consider the singular rank one perturbation $
-\Delta+b<\delta,\cdot>\delta(x),$ where $A_0=-\Delta$
$(\mathcal{D}(A_0)=W_2^2(\mathbb{R}^n)$ is the Laplace operator in
${\mathfrak H}=L_2(\mathbb{R}^n)$ and the associated symmetric
operator $
 A_{\mathrm{sym}}=-\Delta\upharpoonright{\{\,u(x)\in{W_2^{2}(\mathbb{R}^n)}
  :\, u(0)=0 \,\}}.$

The operator $A_0$ is $t^{-2}$-homogeneous with respect to the set
of scaling transformations ${\mathfrak U}=\{U_t\}_{t\in{(0,
\infty)}}$ in $L_2(\mathbb{R}^n)$, where $U_tf(x)={t}^{n/2}f(tx)$.
Furthermore,  the singular element $\psi=\delta$ is
${t}^{-n/2}$-invariant (cf. \cite{AL1}).

If $n=1$, then $\delta(x)\in{\mathfrak
H}_{-1}(A_0)=W_2^{-1}(\mathbb{R})$, the equation (\ref{ada7}) has a
unique solution and by Proposition~\ref{pip1} the free Laplace
operator $-\Delta$ coincides with the Krein-von Neumann extension
$A_N$ of $A_{\mathrm{sym}}$. The Friedrichs extension $A_F$ has the
form
${A}_F=-d^2/dx^2\upharpoonright{\{\,u(x)\in{W_2^{2}(\mathbb{R}\setminus\{0\})\cap{W_2^{1}(\mathbb{R})}}
 :\, u(0)=0 \,\}}.$

If $n=2$, then (\ref{ada7}) has no solutions and there exists the
unique nonnegative self-adjoint extension $-\Delta=A_N=A_F$ of
$A_{\mathrm{sym}}$.

If $n=3$, then
$\delta(x)\in{W_2^{-2}(\mathbb{R}^{3})}\setminus{W_2^{-1}(\mathbb{R}^{3})}$,
the equation (\ref{ada7}) has a unique solution and $-\Delta=A_F$.
The Krein-von Neumann extension ${A_N}$ has the form
$$
{A}_Nf(x)=-\Delta{u}(x)-{u(0)}\frac{e^{-|x|}}{|x|}, \
\mathcal{D}(A_N)=\{\,f=u(x)+{u(0)}\frac{e^{-|x|}}{|x|} :\,
 {u}\in{W_2^2}({\mathbb R}^3)\}.
$$

Another description of the Krein-von Neumann extension of
$A_{\mathrm{sym}}$ obtained with the aid of the Fourier
transformation can be founded in \cite{AT}.
\end{example}

\section{Operator realizations in the case of singular perturbations with symmetries}
\label{sec5}

In this section, operator realizations $A_{\mathbf{B}}$ of
(\ref{ne3}) given by formula (\ref{k41}) are studied under the
condition that the unperturbed operator $A_0$ and the singular
elements $\psi_j$ in (\ref{ne3}) are, respectively,
$p(t)$-homogeneous and $\xi_j(t)$-invariant with respect to
${\mathfrak U}$.

\subsection{$p(t)$-Homogeneous operator realizations.}

\begin{theorem}\label{p5}
Let an admissible operator $\widetilde{A}$ for the regularization of
(\ref{ne3}) be chosen to be $p(t)$-homogeneous. Then the operator
$A_{\mathbf{B}}$ defined by (\ref{k41}) is $p(t)$-homogeneous if and
only if the relations
\[
\xi_i(t)\xi_j(t)=p(t), \qquad \forall{t}\in{\mathfrak T}
\]
hold for all indices ${1}\leq{i,j}\leq{n}$ corresponding to non-zero
entries $b_{ij}$ of ${\mathbf{B}}$.
\end{theorem}

\begin{proof} By Lemma \ref{l56}, the operator $A_{\mathrm{sym}}^*$
is $p(t)$-homogeneous. Hence, in view of (\ref{les102}),
$A_{\mathbf{B}}$ is $p(t)$-homogeneous if and only if $U_{g(t)} :
{\mathcal D}(A_{\mathbf{B}}) \to {\mathcal D}(A_{\mathbf{B}}), \quad
 \forall{t}\in\mathfrak{T}.$  By (\ref{k41}), this relation can be
 rewritten as
 \begin{equation}\label{sas99}
 {\mathbf{B}}\Gamma_0U_{g(t)}f=\Gamma_1U_{g(t)}f, \ \ \ \
 \forall{t}\in\mathfrak{T}, \ \ \ \forall{f}\in\mathcal{D}(A_{\mathbf{B}}).
 \end{equation}

 Since the admissible operator $\widetilde{A}$ is $p(t)$-homogeneous,
 the boundary operator $\Gamma_0$ satisfies (\ref{les101}). Therefore,
 ${\mathbf{B}}\Gamma_0U_{g(t)}f={\mathbf{B}}{\mathbf{\Xi}}(t)\Gamma_0f$.
 On the other hand, relations (\ref{lesia101}) and (\ref{les10})
 lead to the equality
\begin{equation}\label{sas102}
{\Gamma}_1U_{g(t)}f=p(t){\mathbf{\Xi}}^{-1}(t){\Gamma}_1f, \quad
\forall{f\in{\mathcal D}(A_{\mathrm{sym}}^*)}.
\end{equation}
The last two equalities and (\ref{k41}) show that the relation
(\ref{sas99}) is equivalent to the matrix equality
${\mathbf{\Xi}}(t){\mathbf{B}}{\mathbf{\Xi}}(t)=p(t){\mathbf{B}}$,
${t}\in\mathfrak{T}.$ Rewriting this componentwise, one
 obtains the equalities
 $\xi_i(t)\xi_j(t)b_{ij}=p(t)b_{ij}, \ {{1}\leq{i,j}\leq{n}}$.
 \end{proof}
 \begin{cor}\label{c45}
 If there exists a point $t_0\in{\mathfrak T}$ such that
 $p(t_0)\not=1$ and relations $\xi_i(t_0)\xi_j(t_0)=p(t_0)$ hold
 for all indices ${1}\leq{i,j}\leq{n}$ corresponding to
 non-zero entries $b_{ij}$ of $\mathbf{B}$, then: (i)
 the point $\lambda=0$ belongs to the essential spectrum of $A_{\mathbf{B}}$
and $ \lambda\in\sigma(A_{\mathbf{B}}) \iff
\lambda{p(t_0)}^n\in\sigma(A_{\mathbf{B}}), \ n\in\mathbb{Z};$ \
(ii) \ the operator $A_{\mathbf{B}}$ is nonnegative if and only if
the matrix ${\mathbf B}$ is Hermitian.
 \end{cor}
\begin{proof}
If the matrix $\mathbf{B}$ satisfies the conditions above, then
 $A_{\mathbf{B}}$ is $p(t)$-homogeneous with respect to the family
${\mathfrak U}_0:=\{\,U_t\in {\mathfrak U}:\, t\in\{t_0,
g(t_0)\}\,\}$. Now, to establish (i), it suffices to use Lemma
\ref{prop33} with $A=A_{\mathbf{B}}$.

Obviously, the matrix ${\mathbf B}$ is Hermitian if and only if the
operator $A_{\mathbf{B}}$ defined by (\ref{k41}) is self-adjoint.
Using Lemma \ref{t33} and Theorem \ref{p5} one derives (ii).
\end{proof}

\begin{proposition}\label{new2008}
Assume that the  singular elements $\psi_j$ in (\ref{ne3}) form a
${\mathfrak H}_{-1}(A_0)$-independent orthonormal system in
${\mathfrak H}_{-2}(A_0)$, the system (\ref{kak9}) has a unique
solution $\mathbf{R}$, and a $p(t)$-homogeneous admissible operator
$\widetilde{A}$ is chosen for the regularization of (\ref{ne3}).
Then a self-adjoint operator realization $A_{\mathbf{B}}$ of
(\ref{ne3}) is nonnegative if and only if \ $
\det({\mathbf{BR}}+\mathbf{E})\not=0 \quad \mbox{and} \quad
0\leq-({\mathbf{BR}}+\mathbf{E})^{-1}{\mathbf{B}}\leq-\mathbf{R}^{-1},
$ \ where $\mathbf{E}$ stands for the identity matrix.
\end{proposition}
\begin{proof}
By Theorem \ref{new2007}, the Krein-von Neumann extension $A_N$ of
$A_{\mathrm{sym}}$ coincides with a $p(t)$-homogeneous admissible
operator $\widetilde{A}$ and it is defined by (\ref{lesia200}),
where $\mathbf{R}$ is the solution of (\ref{kak9}). Furthermore, the
Friedrichs extension $A_F$ coincides with $A_0$. Combining these
observations with \cite[Theorem 3]{wife} the statement follows. For
completeness some of the details are repeated here.

By (\ref{sas44}) a self-adjoint operator $A_{\mathbf{B}}$ is
nonnegative if and only if $-1\in\rho(A_{\mathbf{B}})$ and
\begin{equation}\label{as55}
0\leq{C_{\mathbf{B}}}\leq{C_N},
\end{equation}
where $C_{\mathbf{B}}=(A_{\mathbf{B}}+I)^{-1}-(A_0+I)^{-1}$ and
$C_{N}=(A_{N}+I)^{-1}-(A_0+I)^{-1}$ are self-adjoint operators in
$\mathcal{H}=\ker(A_{\mathrm{sym}}^*+I)$.

It follows from (\ref{lesia101}) and (\ref{k41}) that
\begin{equation}\label{as7}
\mathcal{D}(A_{\mathbf{B}})=\{\, f\in\mathcal{D}(A_{\mathrm{sym}}^*)
:\,
 {\mathbf B}\widehat{\Gamma}_1f=-({\mathbf{BR}}+\mathbf{E})\widehat{\Gamma}_0f
 \,\}.
\end{equation}

Relations (\ref{k9}) and (\ref{as7}) imply
$-1\in\rho(A_{\mathbf{B}})\iff\mathcal{D}(A_{\mathbf{B}})\cap\mathcal{H}=\{0\}\iff\det({\mathbf{BR}}+\mathbf{E})\not=0$.
Since the elements $\psi_j$ are orthonormal in ${\mathfrak
H}_{-2}(A_0)$, the corresponding vectors $h_j$ in (\ref{kk41}) form
an orthonormal basis of $\mathcal{H}$. In that case, the domain
$\mathcal{D}(A_{\mathbf{B}})$ can be also presented as
$\mathcal{D}(A_{\mathbf{B}})=\{\,
f\in\mathcal{D}(-\Delta_{\mathrm{sym}}^*) :\,
 \mathbf{C}_{\mathbf B}\widehat{\Gamma}_1f=\widehat{\Gamma}_0f
 \,\}$, where $\mathbf{C}_{\mathbf B}$ is the matrix representation
 of $C_{\mathbf{B}}$ with respect to the basis $\{h_j\}_1^n$.
Comparing this with (\ref{as7}) one gets  $\mathbf{C}_{\mathbf
B}=-({\mathbf{BR}}+\mathbf{E})^{-1}{\mathbf{B}}$.

Similar reasonings for the operator $A_N$ defined by
(\ref{lesia200}) give $\det{\mathbf{R}}\not=0$ (since
$-1\in\rho(A_{N})$) and $\mathbf{C}_N=-\mathbf{R}^{-1}$ . By
substituting the obtained expressions for $\mathbf{C}_{\mathbf B}$
and $\mathbf{C}_N$ into (\ref{as55}) one completes the proof.
\end{proof}

\begin{remark} A description of nonnegative self-adjoint operator
realizations of (\ref{ne3}) given above is based on the specific
form of boundary operators $\Gamma_i$. A general approach to the
description of nonnegative self-adjoint extensions of a symmetric
operator has been proposed recently in \cite{AT}.
\end{remark}

\subsection{The Weyl function and the resolvent formula.}
Let $({\mathbb C}^n, \Gamma_0, \Gamma_1)$ be the boundary triplet of
$A_{\mathrm{sym}}^*$ constructed in Lemma \ref{l23} and let
$\widetilde{A}$ be a self-adjoint extension of $A_{\mathrm{sym}}$
defined by (\ref{lesia200}).

The $\gamma$-field $\gamma(z)$ and the Weyl function $\mathbf{M}(z)$
associated with the boundary triplet $({\mathbb C}^n, \Gamma_0,
\Gamma_1)$ are defined by
\begin{equation}\label{sas96}
\gamma(z)=(\Gamma_0\upharpoonright {{{\mathcal H}}_z})^{-1}, \qquad
\mathbf{M}(z)=\Gamma_1\gamma(z), \quad z\in\rho(\widetilde{A}).
\end{equation}
Here ${\mathcal H}_z=\ker(A_{\mathrm{sym}}^*-zI)$, $z\in\mathbb{C}$
denote the defect subspaces of $A_{\mathrm{sym}}$. The mappings
$\Gamma_i$ are defined by (\ref{k9}) and $\mathbf{M}(z)$ is an
$n\times{n}$-matrix function.

\begin{theorem}\label{pop1}
The operator $\widetilde{A}$ is $p(t)$-homogeneous if and only if
for at least one point ${z=z_0}\in{\mathbb C}\setminus{\mathbb R}$
(and then for all non-real points $z$) the Weyl function
$\mathbf{M}(z)$ satisfies the relation
\begin{equation}\label{les3}
 p(t)\mathbf{M}(z)={\mathbf{\Xi}}(t)\mathbf{M}(p(t)z){\mathbf{\Xi}}(t),
 \quad  \forall{t}\in{\mathfrak T},
\end{equation}
where ${\mathbf{\Xi}}(t)$ is defined by (\ref{did4}).
\end{theorem}

\begin{proof}
Let $f_z\in{\mathcal H}_z$, $z\in\mathbb{C}$. Then Lemma
\ref{prop33} and relation (\ref{tat1}) imply
\begin{equation}\label{sas101}
U_{g(t)}f_z\in\ker(A_{\mathrm{sym}}^*-\frac{z}{p(g(t))}I)=\ker(A_{\mathrm{sym}}^*-p(t)zI)={\mathcal{H}_{p(t)z}}.
\end{equation}

Putting $f=f_z\in{\mathcal H}_z$ in (\ref{sas102}), using
(\ref{sas101}), and observing that
$\mathbf{M}(z)\Gamma_0f_z=\Gamma_1f_z$, $z\in{\mathbb{C}}$ (see
(\ref{sas96})), one can rewrite (\ref{sas102}) as follows:
\begin{equation}\label{sas103}
\mathbf{M}(p(t)z)\Gamma_0U_{g(t)}f_z=p(t){\mathbf{\Xi}}^{-1}(t)\mathbf{M}(z)\Gamma_0f_z.
\end{equation}

If the identity (\ref{les3}) holds for some non-real $z=z_0$, then
(\ref{sas103}) implies that
\begin{equation}\label{sas104}
\Gamma_0U_{g(t)}f=\mathbf{\Xi}(t)\Gamma_0f
\end{equation}
for all $f=f_{z_0}\in{\mathcal H}_{z_0}$. Since
$\mathbf{M^*}(z)=\mathbf{M}(\overline{z})$ \cite{DM} and hence,
(\ref{les3}) holds for $\overline{z}_0$, the  relation
(\ref{sas104}) is also true for $f=f_{\overline{z}_0}\in{\mathcal
H}_{\overline{z}_0}$. Moreover, (\ref{sas104}) holds for all
$f\in\mathcal{D}(A_{\mathrm{sym}})$ since
$\Gamma_0f=\Gamma_0U_{g(t)}f=0$ by (\ref{e7}). Consequently,
(\ref{sas104}) is true on the domain ${\mathcal
D}(A_{\mathrm{sym}}^*)=\mathcal{D}(A_{\mathrm{sym}})\dot{+}{\mathcal
H}_{z_0}\dot{+}{\mathcal H}_{\overline{z}_0}$. By Theorem \ref{t5}
this provides the $p(t)$-homogeneity of $\widetilde{A}$.

Conversely, assume that $\widetilde{A}$ is $p(t)$-homogeneous. In
this case, (\ref{sas104}) holds for all $f\in{\mathcal
D}(A_{\mathrm{sym}}^*)$ (see (\ref{les101})). But then, for all
non-real $z$ and all $f_z\in{\mathcal{H}_z}$,
\begin{eqnarray*}
\mathbf{M}(p(t)z)\mathbf{\Xi}(t)\Gamma_0f_z\stackrel{(\ref{sas104})}{=}\mathbf{M}(p(t)z)\Gamma_0U_{g(t)}f_z\stackrel{(\ref{sas101})}{=}
\Gamma_1U_{g(t)}f_z & &  \\
\stackrel{(\ref{sas102})}{=}p(t){\mathbf{\Xi}}^{-1}(t){\Gamma}_1f_z{=}
p(t){\mathbf{\Xi}}^{-1}(t)\mathbf{M}(z)\Gamma_0f_z
\end{eqnarray*}
that justifies (\ref{les3}). Theorem  \ref{pop1} is proved.
\end{proof}

Let $A_{\mathbf{B}}$ be a self-adjoint realization of (\ref{ne3})
defined by (\ref{k41}). Then the resolvents of $A_{\mathbf{B}}$ and
$\wt A$ are connected via Krein's formula
\begin{equation}\label{les34}
(A_{\mathbf{B}}-zI)^{-1}=(\widetilde{A}-zI)^{-1}+\gamma(z)({\mathbf{B}}-{\mathbf{M}}(z))^{-1}\gamma(\overline{z})^{*},
\quad z\in\rho(A_{\mathbf{B}})\cap\rho(\wt A).
\end{equation}

The explicit form of $\mathbf{M}(z)$ can be found as follows. By
(\ref{lesia101}) it is easy to see that the Weyl functions
$\mathbf{M}(z)$ and ${\mathbf{\wh M}}(z)$ associated with the
boundary triplets (\ref{k9}) and (\ref{lesia99}), respectively, are
connected via the linear fractional transform
\begin{equation}\label{ura1}
  \mathbf{M}(z)=-(\mathbf{R}+\mathbf{\wh M}(z))^{-1}, \quad z\in\cmr.
\end{equation}
The boundary triplet (\ref{lesia99}) is one of the most used
boundary triplets and the corresponding Weyl function $\mathbf{\wh
M}(z)$ is studied well. In particular, if the singular elements
$\psi_j$ in (\ref{ne3}) form an orthonormal system in
$\mathfrak{H}_{-2}$, then (see \cite[Remark 4]{DM})
$$
\mathbf{\wh
M}(z)=(z+1)P_{\mathcal{H}}[I+(z+1)(A_0-zI)^{-1}]P_{\mathcal{H}}.
$$
By combining this relation with (\ref{ura1}) one gets an explicit
form for $\mathbf{M}(z)$.
\begin{example}\label{Examp}
{\it A point interaction for \textsf{p}-adic Schr\"{o}dinger type
operator.}
 Let  \textsf{p} be a fixed prime number and let
$\mathbb{Q}_\textsf{p}$ be the field of \textsf{p}-adic numbers. The
operation of differentiation is not defined in the \textsf{p}-adic
analysis of complex-valued functions defined on
$\mathbb{Q}_\textsf{p}$ and the Vladimirov operator of the
fractional \textsf{p}-adic differentiation
$$
D^{\alpha}f(x)=\frac{\textsf{p}^\alpha-1}{1-\textsf{p}^{-1-\alpha}}\int_{\mathbb{Q}_\textsf{p}}\frac{f(x)-f(y)}{|x-y|_\textsf{p}^{1+\alpha}}d\mu(y),
\quad \alpha>0
$$
is used as an analog of it (see \cite{KO1} for details). Here
$|\cdot|_\textsf{p}$ and $d\mu(y)$ are, respectively, the
\textsf{p}-adic norm and the Haar measure on
$\mathbb{Q}_\textsf{p}$. The operator $D^{\alpha}$ is positive and
self-adjoint in the Hilbert space $L_2(\mathbb{Q}_\textsf{p})$ of
complex-valued square integrable functions on
$\mathbb{Q}_\textsf{p}$.
 $\textsf{P}$-adic Schr\"{o}dinger-type operators with potentials
$V(x) : \mathbb{Q}_\textsf{p}\to\mathbb{C}$ are defined as
$D^{\alpha}+V(x)$.

Denote $\mathfrak{T}=\{t=\textsf{p}^n : n\in\mathbb{Z}\}$ and
consider a family $\mathfrak{U}=\{U_t\}_{t\in\mathfrak{T}}$ of
unitary operators $U_tf(x)={t}^{-1/2}f(tx)$ acting in
$L_2(\mathbb{Q}_\textsf{p})$. Obviously, $U_t$ satisfies (\ref{e1})
with the function of conjugation $g(t)=1/t$, c.f. (\ref{ee7}). It
follows from \cite{KOZ} that $U_tD^{\alpha}=t^{\alpha}D^{\alpha}U_t,
 t\in\mathfrak{T}$. Hence, $D^{\alpha}$ is $t^\alpha$-homogeneous
with respect to $\mathfrak{U}$.

Since $D^\alpha$ is a $\textsf{p}$-adic pseudo-differential operator
its domain of definition $\mathcal{D}(D^{\alpha})$ need not contain
functions continuous on $\mathbb{Q}_\textsf{p}$ and, in general, may
happen that the formal expression
\begin{equation}\label{as1}
D^{\alpha}+b<\delta, \cdot>\delta(x), \qquad  b\in\mathbb{R}
\end{equation}
and the associated symmetric operator
$A_{\mathrm{sym}}=D^\alpha\upharpoonright{\{\,u(x)\in\mathcal{D}(D^{\alpha})
  :\, u(0)=0 \,\}}$
are not defined on $\mathcal{D}(D^{\alpha})$. It is known \cite{KT}
that the domain $\mathcal{D}(D^{\alpha})$ consists of continuous
functions on $\mathbb{Q}_\textsf{p}$ and the Dirac delta function
$\delta(x)$ is well-defined on
$\mathfrak{H}_2(D^\alpha)=\mathcal{D}(D^{\alpha})$ if and only if
$\alpha>1/2$. Furthermore, $\delta(x)$ is $\sqrt{t}$-invariant with
respect to $\mathfrak{U}$ and
$\delta(x)\in\mathfrak{H}_{-2}(D^\alpha)\setminus\mathfrak{H}_{-1}(D^\alpha)$
if $1/2<\alpha\leq{1}$, while
$\delta(x)\in\mathfrak{H}_{-1}(D^\alpha)$ if $\alpha>1$.

It follows from \cite[Lemma 3.7]{KO1} and \cite[Lemma 2.1]{KT} that
$$
h(x)=(D^\alpha+I)^{-1}\delta=\sum_{N=-\infty}^{\infty}\sum_{j=1}^{\textsf{p}-1}\textsf{p}^{-N/2}\big[\textsf{p}^{\alpha(1-N)}+1\big]^{-1}\psi_{Nj0}(x),
$$
where the functions $\psi_{Nj0}(x)$ ($N\in\mathbb{Z},
j=1\ldots,\textsf{p}-1$) form a part of the $\textsf{p}$-adic
wavelet basis $\{\psi_{Nj\epsilon}(x)\}$ recently constructed in
\cite{KOZ}.

The equation (\ref{ada7}) takes the form
\begin{equation}\label{ada8}
(t-t^\alpha)r=\sqrt{t}(1-t^\alpha))(h, U_th), \qquad
\forall{t}\in\mathfrak{T}.
\end{equation}

A simple analysis shows that (\ref{ada8}) has no solutions for
$\alpha=1$. In that case the initial operator $D^1$ is a unique
nonnegative self-adjoint extension of $A_{\mathrm{sym}}$, see
Proposition \ref{pip1}. If  $\alpha\not=1$ ($\alpha>1/2$), then
(\ref{ada8}) has a unique solution $r\in\mathbb{R}$ that determines
a unique $t^\alpha$-homogeneous admissible operator $\widetilde{A}$
for the regularization of (\ref{as1}) by the formula (cf.
(\ref{lesia200}))
$$
\widetilde{A}f(x)=D^{\alpha}u(x)+\frac{u(0)}{r}h(x), \quad
\mathcal{D}(\widetilde{A})=\{\, f=u(x)-\frac{u(0)}{r}h(x) :\,
{u}\in\mathcal{D}(D^\alpha)\}.
$$

In view of Proposition \ref{pip1}, the operator $\widetilde{A}$
coincides with the Krein-von Neumann  (Friedrichs) extension of
$A_{\mathrm{sym}}$ for $1/2<\alpha<1$ (resp. for $\alpha>1$).

Let $({\mathbb C}^n, \Gamma_0, \Gamma_1)$ be the boundary triplet of
$A_{\mathrm{sym}}^*$ constructed in Lemma \ref{l23} so that
$\ker\Gamma_0=\mathcal{D}(\widetilde{A})$. By Theorem \ref{ttt12},
self-adjoint operator realizations of (\ref{as1}) in
$L_2(\mathbb{Q}_\textsf{p})$ have the form
$A_bf=A_b(u+{c}h)=D^{\alpha}u-ch, \quad
\forall{u}\in\mathcal{D}(D^{\alpha}),$ where the parameter $c=c(u,
b)\in\mathbb{C}$ is uniquely determined by the relation
$bu(0)=-c[1+br]$. Since $\xi^2(t)=t\neq{t^\alpha}=p(t)$
($\alpha\neq{1}$), Theorem \ref{p5} shows that $A_b$ is
$t^\alpha$-homogeneous if and only if $b=0$ or $b=\infty$
($A_{\infty}\equiv\widetilde{A}$).

Let $\alpha>1$. It follows from \cite{AKT} that the Weyl function
associated with $({\mathbb C}^n, \Gamma_0, \Gamma_1)$ has the form
$$
\mathbf{M}(z)=-\frac{1}{(\textsf{p}-1)\sum_{N=-\infty}^{\infty}\frac{\textsf{p}^{-N}}{\textsf{p}^{\alpha(1-
N)}-z}}.
$$

By virtue of Theorem \ref{pop1}, $\mathbf{M}(z)$ satisfies the
relation $t^{\alpha-1}\mathbf{M}(z)=\mathbf{M}(t^\alpha{z})$,
$\forall{t}\in\mathfrak{T}$. This simplifies the spectral analysis
of $A_b$, see \cite{AKT} for details.
\end{example}

\section{Schr\"{o}dinger operators with singular perturbations $\xi(t)$-invariant
 with respect to scaling transformations in ${\mathbb R^3}$}
 \label{sec6}
It is well known (see, e.g. \cite{AL1, CF}) that the Schr\"{o}dinger
operator $A_0=-\Delta$, \ $(D(\Delta)=W_2^2(\mathbb{R}^3))$  is
$t^{-2}$-homogeneous with respect to the set of scaling
transformations ${\mathfrak U}=\{U_t\}_{t\in{(0, \infty)}}$
($U_tf(x)={t}^{3/2}f(tx)$) in $L_2(\mathbb{R}^3)$. It is clear that
$U_t$ satisfies (\ref{e1}) with the function of conjugation
$g(t)=1/t$.

The elements $U_t$ of ${\mathfrak U}$ possess the additional
multiplicative property $U_{t_1}U_{t_2}=U_{t_2}U_{t_1}=U_{t_1t_2}$
that enables one to describe all measurable functions $\xi(t)$ for
which there exist $\xi(t)$-invariant singular elements
$\psi\in{W}_2^{-2}(\mathbb{R}^3)$.

\begin{theorem}\label{dod1}
Let $\xi(t)$ be a real measurable function defined on $(0,\infty)$.
Then $\xi(t)$-invariant singular elements
$\psi\in{W}_2^{-2}(\mathbb{R}^3)\setminus{L_2(\mathbb{R}^3)}$ exist
if and only if $\xi(t)=t^{-\alpha}$, where $0<\alpha<2$.
\end{theorem}

\begin{proof} Let $\psi\in{W}_2^{-2}(\mathbb{R}^3)\setminus{L_2(\mathbb{R}^3)}$
be $\xi(t)$-invariant with respect to ${\mathfrak U}$. Since
$U_{t_1}U_{t_2}=U_{t_2}U_{t_1}=U_{t_1t_2}$, equality (\ref{tat10})
gives $\xi(t_1)\xi(t_2)=\xi(t_1t_2)$ $(t_i>0)$ that is possible only
if $\xi(t)=0$ or $\xi(t)=t^{-\alpha}$ $(\alpha\in\mathbb{R})$
\cite[Chap.IV]{HP}. Furthermore, Proposition \ref{p1} enables one to
restrict the set of possible functions $\xi(t)$ as follows:
$\xi(t)=t^{-\alpha}$, where $0<\alpha<2$.

To complete the proof of Theorem \ref{dod1} it suffices to construct
$t^{-\alpha}$-invariant singular elements for $0<\alpha<2$.

Fix $m(w)\in{L_2(S^2)}$, where $L_2(S^2)$ is the Hilbert space of
square-integrable functions on the unit sphere $S^2$ in
$\mathbb{R}^3$, and determine the functional
$\psi(m,\alpha)\in{W}_2^{-2}({\mathbb R}^3)$ by the formula
\begin{equation}\label{ada66}
 <\psi(m,\alpha), u>=
 \int_{{\mathbb R}^3}\frac{m(w)}{|y|^{3/2-\alpha}(|y|^2+1)}(|y|^2+1){\widehat{u}}(y)dy
  \quad (y=|y|w\in\mathbb{R}^3),
\end{equation}
where ${\widehat{u}}(y)=\frac{1}{(2\pi)^{3/2}}\int_{{\mathbb
R}^3}e^{ix\cdot{y}}u(x)dx$ is the Fourier transformation of
$u(\cdot)\in{W}_2^{2}({\mathbb R}^3)$.

It is easy to verify that
\begin{equation}\label{kaka2}
{\widehat{(U_{g(t)}u)}}(y)=
\widehat{(U_{1/t}u)}(y)=\frac{1}{(2\pi{t})^{3/2}}\int_{{\mathbb
 R}^3}e^{iy\cdot{x}}u(x/t)dx=U_{t}{\widehat{u}}(y)=t^{3/2}{\widehat{u}}(ty).
\end{equation}
Using (\ref{ada66}) and (\ref{kaka2}), one obtains \ $
<\psi(m,\alpha), U_{g(t)}u>=t^{-\alpha}<\psi(m,\alpha), u>$ for all
${u}\in{{W}_2^{2}({\mathbb R}^3)}$.  By (\ref{tat3}) this means that
the functional $\psi(m,\alpha)$ is $t^{-\alpha}$-invariant with
respect to ${\mathfrak U}$. Theorem \ref{dod1} is proved.
\end{proof}

A more detailed study of functionals that are
$t^{-\alpha}$-invariant with respect to scaling transformations and
the results of \cite{SW} lead to the conclusion that the collection
$\mathcal{L}_{\alpha}$ of all $t^{-\alpha}$-invariant singular
elements
$\psi\in{W}_2^{-2}(\mathbb{R}^3)\setminus{L_2(\mathbb{R}^3)}$ can be
described as follows: \ $
 \mathcal{L}_{\alpha}=\left\{\psi=\psi(m,\alpha)\,
 :\, {m(w)}\in{L}_2(S^2), \quad m(w)\not=0 \right\}.$

Let us consider the formal expression
\begin{equation}\label{eee3}
-\Delta+\sum_{i,j=1}^{n}{b}_{ij}<\psi_j,\cdot>\psi_i, \qquad
 b_{ij}\in\mathbb{C}, \quad  n\in\mathbb{N},
 \end{equation}
where all singular elements $\psi_j$ are assumed to be
$t^{-\alpha}$-invariant with respect to scaling transformations for
 a fixed  $\alpha$, i.e.,  $\psi_j=\psi(m_j,\alpha)$.
The symmetric operator $A_{\mathrm{sym}}=-\Delta_{\mathrm{sym}}$
associated with (\ref{eee3}) takes the form
\begin{equation}\label{tato72}
 -\Delta_{\mathrm{sym}}=-\Delta\upharpoonright_{{\mathcal{D}}(\Delta_{\mathrm{sym}})},
   \ {\mathcal{D}}(\Delta_{\mathrm{sym}})=\{\,u(x)\in{W_2^2(\mathbb{R}^3)}  :\,  <\psi_j, u>=0, \
 1\leq{j}\leq{n}\,\},
\end{equation}
where $<\psi_j, u>$ are defined by (\ref{ada66}).

Comparing (\ref{ada5}) and (\ref{ada66}), one sees that the
functions $h_j=(\mathbb{A}_0+I)^{-1}\psi(m_j,\alpha)$ in
(\ref{kk41}) have the form
\begin{equation}\label{kaka1}
h_j(x)=\left(\frac{\overline{m_j(w)}}{|y|^{3/2-\alpha}(|y|^2+1)}\right)^{\displaystyle\stackrel{\lor}{}}(x)=\overline{\left(\frac{m_j(w)}{|y|^{3/2-\alpha}(|y|^2+1)}\right)^{\displaystyle\stackrel{\land}{}}(x)},
\end{equation}
where the symbol $\stackrel{\lor}{}$ denotes the inverse Fourier
transformation.

A simple analysis of (\ref{kaka1}) shows that
$h_j\in{L_2(\mathbb{R}^3)}\setminus{W_2^1(\mathbb{R}^3)}$ for
$1\leq\alpha<2$ and $h_j\in{W_2^1(\mathbb{R}^3)}$ for $0<\alpha<1$.
In the latter case, Corollary \ref{p121} and Lemma \ref{l34} imply
that the Friedrichs extension $-\Delta_F$ is a unique
$t^{-2}$-homogeneous admissible operator for the regularization of
(\ref{eee3}).

\begin{proposition}\label{t66}
Let $1<\alpha<2$. Then the Krein-von Neumann extension $-\Delta_N$
of $-\Delta_{\mathrm{sym}}$ is a unique $t^{-2}$-homogeneous
admissible operator for the regularization of (\ref{eee3}).
\end{proposition}
\begin{proof}
If $1<\alpha<2$, then all the elements $\psi_j$ in \eqref{eee3} are
$W_2^{-1}(\mathbb{R}^3)$-independent. Let us show that the system
(\ref{kak9}) has a unique solution $\mathbf{R}=(r_{ij})_{i,j=1}^n$
that does not depend on $t>0$. Since the both parts of (\ref{kak9})
are equal to zero for $t=1$, one can suppose that $t>0$ and
$t\not=1$.

It follows from (\ref{kaka2}) and (\ref{kaka1}) that
\[
\begin{split}
\overline{U_th_i(x)}
 &=U_t\left(\frac{m_i(w)}{|y|^{3/2-\alpha}(|y|^2+1)}\right)^{\displaystyle\stackrel{\land}{}}(x)
 =\left(U_{1/t}\frac{m_i(w)}{|y|^{3/2-\alpha}(|y|^2+1)}\right)^{\displaystyle\stackrel{\land}{}}(x)  \\
 &=t^{2-\alpha}\left(\frac{m_i(w)}{|y|^{3/2-\alpha}(|y|^2+t^2)}\right)^{\displaystyle\stackrel{\land}{}}(x).
\end{split}
\]

Hence,
\[
\begin{split}
 (h_j,U_th_i)
 &=t^{2-\alpha}\int_{\mathbb{R}^3}\frac{m_i(w)\overline{m_j(w)}}{|y|^{3-2\alpha}(|y|^2+t^2)(|y|^2+1)}\,dy \\
&
=(m_i,m_j)_{L_2}\int_{0}^{\infty}\frac{t^{2-\alpha}}{|y|^{1-2\alpha}(|y|^2+t^2)(|y|^2+1)}d|y|
\\ &  =c_\alpha\frac{t^{\alpha}-t^{2-\alpha}}{t^2-1}(m_i,m_j)_{L_2},
\end{split}
\]
where $c_\alpha=\int_0^\infty\frac{|y|^{3-2\alpha}}{|y|^2+1}d|y|$
and $(m_i, m_j)_{L_2}=\int_{S^2}m_i(w)\overline{m_j(w)}dw$ is the
scalar product in $L_2(S^2)$.  Substituting the expression for
$(h_j,U_th_i)$ into (\ref{kak9}) one gets a unique solution
$\mathbf{R}=(r_{ij})_{i,j=1}^n$, where $r_{ij}=-c_\alpha(m_i,
m_j)_{L_2}$. By Theorem \ref{new2007}, the obtained solution
determines a unique $t^{-2}$-homogeneous admissible operator
$\widetilde{A}$ for the regularization of (\ref{eee3}) that
coincides with $-\Delta_N$.
\end{proof}

\begin{remark}
If $\alpha=1$, then (\ref{kak9}) has no solution, there are no
$t^{-2}$-homogeneous admissible operators for (\ref{eee3}), and the
Friedrichs $-\Delta=-\Delta_F$ and the Krein-von Neumann $-\Delta_N$
extensions of $-\Delta_{\mathrm{sym}}$ are not transversal.
\end{remark}

\begin{corollary}\label{popa1}
For a fixed $1<\alpha<2$ assume that $\psi_j=\psi(m_j,\alpha)$ in
(\ref{eee3}) form an orthonormal system in $W_2^{-2}(\mathbb{R}^3)$
and self-adjoint operator realizations
$A_{\mathbf{B}}=-\Delta_{\mathbf{B}}$ of (\ref{eee3}) are defined by
(\ref{k41}) with $\ker\Gamma_0=\mathcal{D}(-\Delta_N)$. Then
$-\Delta_{\mathbf{B}}$ is nonnegative if and only if
$\det(\beta_\alpha{\mathbf{B}}-\mathbf{E})\not=0$ and
$0\leq\beta_\alpha{\mathbf{B}}[\beta_\alpha{\mathbf{B}}-\mathbf{E}]^{-1}\leq\mathbf{E}$,
where
\begin{equation}\label{new555}
\beta_\alpha=\left[\int_0^\infty\frac{|y|^{3-2\alpha}}{|y|^2+1}d|y|\right]\left[\int_{0}^{\infty}\frac{1}{|y|^{1-2\alpha}(|y|^2+1)^2}d|y|\right]^{-1}.
\end{equation}
\end{corollary}
\begin{proof}
Since $\psi(m_j,\alpha)$ are orthonormal in $W_2^{-2}(\mathbb{R}^3)$
the functions $h_j(x)$ determined by (\ref{kaka1}) are orthonormal
in $L_2(\mathbb{R}^3)$. This means that $(m_i, m_j)_{L_2}=0$
($i\not=j$) and $(m_i,
m_j)_{L_2}\int_{0}^{\infty}\frac{1}{|y|^{1-2\alpha}(|y|^2+1)^2}d|y|=1$.
The obtained relations allows one to rewrite the unique solution
$\mathbf{R}=-c_\alpha((m_i, m_j)_{L_2})_{i,j=1}^n$ of (\ref{kak9})
in a more explicit form: $\mathbf{R}=-\beta_\alpha\mathbf{E}$, where
$\beta_\alpha$ is defined by (\ref{new555}). Using Proposition
\ref{new2008} one completes the proof.
\end{proof}

Note that the delta function $\delta(\cdot)$ belongs to
$\mathcal{L}_{3/2}$. For this reason, the expression (\ref{eee3})
where all $\psi_j\in\mathcal{L}_{3/2}$ can be considered as a
generalization of the classical one-point interaction
$-\Delta+b<\delta, \cdot>\delta$. In that case the parameter
$\beta_{\alpha}$ in Corollary \ref{popa1} can be easily calculated:
$\beta_{3/2}=2$.

\begin{theorem}\label{ttt5}
Let $\alpha=3/2$. Then for any self-adjoint operator realization
$A_{\mathbf{B}}=-\Delta_{\mathbf{B}}$ of (\ref{eee3}) defined by
(\ref{k41}), the following statements are true:
\begin{enumerate}

\def\labelenumi{\rm (\roman{enumi})}

\item if $-\Delta_{\mathbf{B}}$ is nonnegative, then the wave
operators
$W_{\pm}=\lim_{t\to\pm\infty}e^{-it\Delta_{\mathbf{B}}}e^{i\Delta{t}}$
exist and are unitary operators in $L_2({\mathbb R}^3)$;

\item if $-\Delta_{\mathbf{B}}$ is nonnegative and the singular elements
$\psi_j=\psi(m_j, 3/2)$ in (\ref{eee3}) form an orthonormal system
in $W_2^{-2}(\mathbb{R}^3)$, then the $S$-matrix
$$
 {\mathbb S}_{(-\Delta_{\mathbf{B}},-\Delta)}=FW_{+}^*W_{-}F^{-1}
$$
($F$ is the Fourier transformation in $L_2(\mathbb{R}^3)$) of the
Schr\"{o}dinger equation $iu_{t}=-\Delta_{\mathbf{B}}u$ coincides
with the boundary value
$\mathbb{S}_{(-\Delta_{\mathbf{B}},-\Delta)}(\delta)$
($\delta\in\mathbb{R}$) of the contractive operator-valued function
\begin{equation}\label{as4}
\mathbb{S}_{(-\Delta_{\mathbf{B}},-\Delta)}(z)
=(\mathbf{E}-2iz\mathbf{B})(\mathbf{E}+2iz\mathbf{B})^{-1}, \quad
z\in\mathbb{C}_+
\end{equation}
analytic in the upper half-plane $\mathbb{C}_+$.
\end{enumerate}
 \end{theorem}

\begin{proof}The statements follow from \cite[Theorem
3.3]{KM} and \cite[Section 4]{wife}.
\end{proof}

\begin{remark} In \cite{wife} the expression (\ref{as4}) was obtained
by using the Lax--Phillips scattering scheme. Another description of
$\mathbb{S}_{(-\Delta_{\mathbf{B}},-\Delta)}(z)$ in terms of the
Krein's resolvent formula was obtained in \cite{AdPa}. In that
paper, the stationary scattering theory approach has been used.
\end{remark}

\noindent \textbf{Acknowledgements.} The authors thank S.~Albeverio,
Yu.~Arlinskii, and L.~Nizhnik for useful discussions. The first
author (S.H.) is grateful for the support from the Research
Institute for Technology of the University of Vaasa. The second
author (S.K.) expresses his gratitude to the Academy of Finland
(projects 208056, 117656) for the support and the Department of
Mathematics and Statistics of the University of Vaasa for the warm
hospitality.

\end{document}